\documentclass[10pt]{article}
\usepackage{amssymb}
 \newtheorem{prop}{Proposition}
 
 \newtheorem{theo}{Theorem} \newtheorem{lemma}{Lemma}
 
 \newtheorem{rem}{Remark}
 \newcommand{\bs}{\left\{}
 \newcommand{\es}{\right.}
 \newcommand{\ba}{\begin{array}}
 \newcommand{\ea}{\end{array}}
 \newcommand{\be}{\begin{equation}}
 \newcommand{\ee}{\end{equation}}

 \def\RR{{\rm I\hspace{-0.50ex}R} }

 \def\CC{\rm \hbox{C\kern-.57em\raise.47ex
 \hbox{$\scriptscriptstyle |$}\kern+0.5 em}}
 
 \title{The linear and non linear Rayleigh-Taylor instability for the quasi-isobaric profile}
 \author{Olivier Lafitte\thanks{Universit{\'e} de Paris XIII, LAGA, 93 430 Villetaneuse}\thanks{CEA/DM2S, Centre
 d'Etudes de Saclay, 91191 Gif sur Yvette Cedex}}
 \date{June 29, 2007}

\begin{document}
 \maketitle

 \begin{abstract}
 We study the stability of the system of the Euler equations in the neighborhood of the stationary solution associated with the quasi isobaric profile in a gravity field. This situation corresponds to a Rayleigh-Taylor type problem with a smooth base
 density profile which goes from $0$ to $\rho_a$ (of Atwood number
 $A=1$) given by the ablation front
 model with a thermal conductivity exponent $\nu>1$.
This linear analysis leads to the study of the Rayleigh equation for the
 perturbation of the velocity at the frequency $k$:
$$-\frac{d}{dx}(\rho_0(x)\frac{d{\underline u}}{dx})+k^2[\rho_0(x)-\frac{g}{\gamma^2}\rho'_0(x)]{\underline u}=0.$$ 
We denote by the terms 'eigenmode and eigenvalue' a $L^2$ solution of
the Rayleigh equation associated with a value of $\gamma$. Let $L_0>0$
be given. The quasi isobaric profile is $\rho_0(x)=\rho_a\xi(\frac{x}{L_0})$, where ${\dot \xi}=\xi^{\nu+1}(1-\xi)$. We prove that there exists $L_m(k)$, such that, for all $0<L_0\leq L_m$, there exists an eigenmode ${\underline u}$ such that the unique associated eigenvalue $\gamma$ is in $[\alpha_1, \alpha_2]$, $\alpha_1>0$. Its limit when $L_0$ goes to zero is $\sqrt{gk}$. We obtain an expansion of $\gamma$ in terms of $L_0$ as follows:
$$\gamma=\frac{\sqrt{gk}}{\sqrt{1+2(\Gamma(1+\frac{1}{\nu}))^{-1}(\frac{2kL_0}{\nu})^{\frac{1}{\nu}}+O((kL_0)^{\mbox{min}(1,
 \frac{2}{\nu})})}}.$$
 We identify in this paper the expression of the next term of the expansion of $\gamma$ in powers of $L_0^{\frac{1}{\nu}}$.\\
Using the existence of a maximum growth rate $\Lambda$ and the
existence of at least one eigenvalue belonging to $]\frac{\Lambda}{2},
\Lambda[$ (thanks to a semiclassical analysis), we perform the
nonlinear analysis of the incompressible Euler system of equations
using the method introduced by Grenier. This generalizes the result of
Guo and Hwang (which was obtained in the case $\rho_0(x)\geq\rho_l>0$)
to the case where $\rho_0\rightarrow 0$ when $x\rightarrow -\infty$
and  $k_0(x)=\frac{\rho'_0(x)}{\rho_0(x)}$ satisfy $k_0$ regular
enough, bounded, and $k_0\rho_0^{-\frac12}$ bounded, which is the case
in the model associated with the quasi-isobaric profile, according to $\nu>\frac12$.

\end{abstract}
\setcounter{section}{-1}
\section{Statement of the
problem and main result} In this paper, we study
a theoretical system of equations deduced of the fluid
dynamics analysis of an ablation front model.
Such models have been studied from a physical
point of view by many authors (see H.J. Kull and S.I. Anisimov \cite{kullA}, V. Goncharov, \cite{Gon},
P. Clavin and L. Masse \cite{Masse}). They can be considered as a generalization
in the ablation case of the Rayleigh-Taylor
instability, studied in the pioneering works of
J.W. Strutt (Lord Rayleigh) \cite{Rayleigh} and
G. Taylor \cite{Taylor}.\\
The Rayleigh equation models the Rayleigh-Taylor
instability. It is obtained by considering the
linearization of the incompressible 2d Euler
equations around the solution  $(\rho_0(x), 0, 0,
p_0(x))$ (density, velocity, pressure) with
$\frac{dp_0}{dx}+\rho_0(x)g=0$.
The system of equations write
\be\label{E}\bs\ba{l}\partial_t\rho + \partial_x(\rho U)+\partial_z(\rho V)=0\cr
\partial_t(\rho U)+\partial_x(\rho U^2+P)+\partial_z(\rho UV)=-\rho g\cr
\partial_t(\rho V)+\partial_x(\rho UV) + \partial_z(\rho V^2+P)=0\cr
\partial_xU+\partial_zV=0\ea\es\ee
Write $\rho =\rho_0+\sigma$, $U=v_1, V=v_2$, $P=p_0+p$, the linearized system is
\be\label{LE}\bs\ba{l}\partial_t\sigma + \frac{d\rho_0}{dx}v_1=0\cr
\rho_0(x)\partial_tv_1+\partial_xp=-
\sigma g\cr
\rho_0(x)\partial_tv_2+\partial_zp=0\cr
\partial_xv_1+\partial_zv_2=0.\ea\es\ee
from which one deduces, using $v_1={\tilde u}e^{ikz}$, the partial differential
equation
$$-\frac{\partial}{\partial x}(\rho_0(x)\frac{\partial}{\partial x}\partial^2_{t^2}{\tilde
u})+k^2\rho_0(x)\partial^2_{t^2}{\tilde
u}= gk^2\rho'_0(x){\tilde u}.$$
Introduce
$$T(x, z, t)=\frac{\rho_0(x)}{\rho(x, z, t)}, Q(x, z, t)=\frac{p(x, z, t)-p_0(x)}{\rho_0(x)}$$
the system (\ref{E}) is equivalent to
\be\label{Emod}\bs\ba{l}\partial_tT+{\vec U}.\nabla T=k_0(x)uT\cr
\partial_t{\vec U}+({\vec u}.\nabla ){\vec U}+T\nabla Q+TQk_0(x){\vec e}_1=(1-T){\vec g}\cr
\mbox{div}{\vec U}=0.\ea\es\ee
It is a consequence of the equality $T\rho_0^{-1}\nabla p=T\nabla Q+k_0TQ{\vec e}_1+T{\vec g}$ and of $\partial_tT+{\vec U}.\nabla T=-\rho^{-2}(\partial_t\rho+{\vec U}\nabla \rho)$.\\
The associated linearized system in the neighborhood of ${\vec U}=0$, $T=1$, $Q=0$ is
$$\bs\ba{l}\partial_tT=k_0(x){\tilde u}\cr
\mbox{div}{\tilde {\vec u}}=0\cr
\partial_t{\tilde {\vec u}}+\rho_0^{-1}\nabla(\rho_0{\tilde Q})+{\tilde T}{\vec g}=0.\ea\es$$
Assume that the perturbation
is written as (the real part of) a normal mode $$e^{\gamma
t}e^{ikz}{\underline u}(x, kL_0),$$where $k$ is
the wavelength of the transversal perturbation
and $\gamma$ is the growth rate in time of this
perturbation. We obtain the Rayleigh equation (\ref{eqrayleigh0}) (see C.
Cherfils, P.A. Raviart and O.L. \cite{CCLARA}):
\be\label{eqrayleigh0}-{d\over
dx}(\rho_0(x){d{\underline u}\over dx})+(k^2\rho_0(x)
-\frac{gk^2}{\gamma^2}\rho'_0(x)){\underline u}(x)=0.\ee
We consider a {\bf family} of density profiles $\rho_0(x)$ such that $\rho_0(x)=\rho^0(\frac{x}{L_0})$, where $L_0$ is a characteristic length of the base solution. In one of
the physical applications, namely the case of the
ICF, its magnitude is $10^{-5}$ meters, hence
allowing us to consider the limit $L_0\rightarrow
0$.\\ We develop here a constructive
method for the study of the modes associated with the Kull-Anisimov
density profile (see B. Helffer and O.L.
\cite{HL2}). The Kull-Anisimov profile $\rho_0$ is given by
\be\rho_0(x)=\rho_a\xi(\frac{x}{L_0}),\ee
where the function $\xi$ is a non constant solution of
\be\label{xi}{\dot \xi}=\xi^{\nu+1}(1-\xi),\ee
$\nu$ is called the
thermal conduction index.\\ Note that this equation
on the density is NOT obtained from the
incompressible Euler equations but from a
compressible model with thermal conduction
introduced by Kull and Anisimov \cite{kullA} and
used for example in \cite{HL2} or in
\cite{Lafittepreprint}. The Kull-Anisimov profile satisfies $\mbox{lim}_{x\rightarrow +\infty}\rho_0(x)=\rho_a$, where $\rho_a$ denotes the density of the
ablated fluid, and the convergence is exponential, whereas
$$\mbox{lim}_{x\rightarrow -\infty}\rho_0(x)=0$$
and the convergence is rational 
($(-x)^{\frac{1}{\nu}}\rho_0(x)\rightarrow C_0>0$ when $x\rightarrow -\infty$). The associated Atwood number is thus $1$. Remark also that all non constant solutions of (\ref{xi}) differ from a translation.\\
This case may be related to
the case of the water waves (the density of air being much smaller
than the density of water). It is thus a limit case in all the
theoretical set-up used for the study of Euler equations for fluids
of different densities.\\
Note that, in this case, the self adjoint operator associated with
the equation
(\ref{eqrayleigh0}) is not coercive in $H^1(\mathbb R)$. The methods of \cite{CCLARA}, \cite{HelLaf} cannot be used directly. Moreover, the properties of $\rho_0$ do not allow us to apply \cite{GH}, because it relies on $\rho_0(x)\geq\rho_l>0$.\\
However, consider $k_0(x)=\frac{\rho'_0(x)}{\rho_0(x)}$ introduced in the abstract. In our case, it is equal to $L_0^{-1}\xi^{\nu}(1-\xi)$, hence it is a continuous bounded function which admits a maximum $L_{eff}^{-1}$, and, for $\nu>\frac12$, $k_0\rho_0^{-\frac12}$ is bounded. These properties are (for a more general profile) what is needed to obtain the nonlinear result.\\
\paragraph{Remarks} Define the function $r(t, \varepsilon)$ through:
\be\label{defifonctionr}\frac{1}{\varepsilon}(\xi(-\frac{t}{\varepsilon}))^{\nu}(1-\xi(-\frac{t}{\varepsilon}))=
\frac{1}{\nu t}+\varepsilon^{\frac{1}{\nu}}t^{-1-\frac{1}{\nu}}r(t,
\varepsilon).\ee There exists $t_0>0$ and $\varepsilon_0>0$ such
that $r(t, \varepsilon)$ is bounded for $t\geq t_0,
0\leq\varepsilon\leq \varepsilon_0$, and has a $C^{\infty}$ expansion in $\varepsilon, \varepsilon^{\frac{1}{\nu}}$. Define $S$ through $$\varepsilon^{\frac{1}{\nu}}S'(t,
\varepsilon)=\varepsilon^{-1}\frac{\xi'(-\frac{t}{\varepsilon})}{\xi(-\frac{t}{\varepsilon})}-\frac{1}{\nu
t}, \quad \mbox{lim}_{t\rightarrow +\infty}S(t, \varepsilon)=0.$$ We
have the identity \be\label{S}\xi(-\frac{t}{\varepsilon})(\frac{\nu
t}{\varepsilon})^{\frac{1}{\nu}}\exp(\varepsilon^{\frac{1}{\nu}}S(t,
\varepsilon))=1\ee which implies that there exists a function $r$
bounded for $t\geq t_0$ and $\varepsilon\leq \varepsilon_0$ such that
$$\exp(-\nu \varepsilon^{\frac{1}{\nu}}S(t, \varepsilon))=1+\varepsilon^{\frac{1}{\nu}}t^{-\frac{1}{\nu}}r(t, \varepsilon).$$\\
Let $u(y)={\underline u}(L_0 y)$. The Rayleigh equation rewrites
\be\label{eqrayleigh}-{d\over dy}(\xi(y){du\over
dy})+(\varepsilon^2\xi(y) -\lambda \varepsilon
\xi'(y))u(y)=0,\ee
 where $\varepsilon
= kL_0$ and $\lambda= \frac{gk}{\gamma^2}$. We will consider this
equation from now on.\\ 
We shall introduce two equivalent versions of this equation, which are:

\begin{enumerate}
\item the
system on $(U_+, V_+)$ such that
$U_+(y, \varepsilon)=u(y, \varepsilon)e^{\varepsilon y}$ and $V_+$ (given by the
first equation of the system below),
$v(y, \varepsilon)=V_+(y\varepsilon)e^{-\varepsilon y}$:
\be\label{infinity}\bs\ba{l}\frac{dU_+}{dy}=\varepsilon
(1-\lambda)U_+ + \frac{\varepsilon}{\xi(y)}V_+\cr
\frac{dV_+}{dy}=\varepsilon ( \lambda+1) V_+ +
\varepsilon (1-\lambda^2)\xi(y)U_+,\ea\es\ee

\item if we introduce $w=\frac{v}{\xi(y)}$, the system on $(u, w)$ is
\be\label{hyp}\bs\ba{l}\frac{du}{dt}=\lambda u - w\cr
\frac{dw}{dt}=(\lambda^2-1)u-\lambda
w+(\frac{1}{\nu t}
+\varepsilon^{\frac{1}{\nu}}S'(t,
\varepsilon))w.\ea\es\ee

\end{enumerate}
The first part of the main result of this
paper was presented in \cite{lafitteX}, and  the case where
$\xi(y)=\xi(1)(y+1)^{-\frac{1}{\nu}}$ for $y\geq 0$ was solved in
\cite{CCLARA}. The case of the global ablation system was treated in \cite{Lafittepreprint} and is published \cite{Indiana}.\\ We finally recall that, if there
exists a solution in $L^2(\mathbb R)$ of (\ref{eqrayleigh}), then
$\lambda$ satisfies the inequality (see \cite{HelLaf})\footnote{It
is a consequence of $\mbox{max}
(\frac{\dot\xi}{\xi})=\frac{\nu^{\nu}}{(\nu+1)^{\nu+1}}$}
\be\label{ineq}\lambda\geq
\mbox{max}(1,\varepsilon\frac{(\nu+1)^{\nu+1}}{\nu^{\nu}}).\ee
The main result of the first part of this paper is
\begin{theo}\label{HG0}
\begin{enumerate}
\item There exists $\varepsilon_0>0$, and $C^0>0$
such that, for all $\varepsilon\in ]0,
\varepsilon_0[$ there exists
$\lambda(\varepsilon)\in [{1\over 2}, {3\over 2}]$
such that the Rayleigh equation
(\ref{eqrayleigh}) admits a bounded solution $u$
for $\lambda=\lambda(\varepsilon)$, which corresponds to the eigenmode
$u$ and the eigenvalue $\gamma(k,
\varepsilon)=\sqrt{\frac{gk}{\lambda(\varepsilon)}}$, and
$\lambda(\varepsilon)$ satisfies

$$|\lambda(\varepsilon) -1|\leq C^0\varepsilon^{1\over \nu}.$$
\item We have the estimate
$$\ba{ll}\lambda(\varepsilon) &=
1+2(\frac{2\varepsilon}{\nu})^{\frac{1}{\nu}}(\Gamma(1+\frac{1}{\nu}))^{-1}+o(\varepsilon^{\frac{1}{\nu}})\cr
&=1+2(\frac{2\varepsilon}{\nu})^{\frac{1}{\nu}}(\Gamma(1+\frac{1}{\nu}))^{-1}+O(\varepsilon^{\alpha})\ea$$
with $\alpha = \mbox{min}(1, \frac{2}{\nu})$.
\item We have the expansion
$$1-\lambda(\varepsilon)=2(\frac{\varepsilon}{\nu})^{\frac{1}{\nu}}(B_0(0))^{-1}[1+2(\frac{\varepsilon}{\nu})^{\frac{1}{\nu}}\frac{C_0(1,
  0)}{(B_0(0))^2}+o(\varepsilon^{\frac{1}{\nu}})],$$
where
$B_0(0)=-2\int_0^{\infty}s^{\frac{1}{\nu}}e^{-2s}ds=-2^{-\frac{1}{\nu}}\Gamma(1+\frac{1}{\nu})$
and $C_0(1, 0)$ is calculated below in Proposition \ref{expan}.
\end{enumerate}
\end{theo}
This result is a result, for $k$ fixed, in the limit $L_0\rightarrow 0$. It writes also, for $k$ fixed
and for $L_0<\frac{\varepsilon_0}{k}$ as
\be\label{gamma} \gamma =
\sqrt{\frac{gk}{1+2(\frac{2\varepsilon}{\nu})^{\frac{1}{\nu}}(\Gamma(1+\frac{1}{\nu}))^{-1}+
o(\varepsilon^{\frac{1}{\nu}})}}.\ee 
Note that, in this case, the order of magnitude of $\gamma -\sqrt{gk}$ is not in $kL_0$ as in \cite{CCLARA}, but the result of \cite{HelLaf}, based on $\rho_0-\rho_a1_{x>0}\in L^{\nu+\theta'}$ for all $\theta'>0$ is pertinent.\\We have also a result for $k$ going to infinity, which can be stated as
\begin{prop}
\label{sc}
a) Any value $\lambda(\varepsilon)$ such that (\ref{eqrayleigh}) has a $L^2$ non zero solution satisfies $\frac{kg}{(\lambda(\varepsilon))^2}\leq \Lambda^2$, where $\Lambda^2=\frac{g}{L_{eff}}$.\\
b) Any sequence $k\rightarrow \frac{\lambda(k)}{k}$ satisfies the following
$$\mbox{lim}_{k\rightarrow +\infty}\frac{\lambda(\varepsilon)}{k}= L_{eff}=\mbox{min}_y\frac{\xi(y)}{\xi'(y)}L_0.$$
\end{prop}
It is proven in  \cite{HelLaf}.\\
Remark that formula (\ref{gamma}) and Proposition \ref{sc} are not in contradiction. They lead to two different stabilizing mechanisms induced by the transition region: one is a low frequency stabilization when $L_0\rightarrow 0$ and the other one is a high frequency stabilizing mechanism when $k\rightarrow +\infty$.
It is important to notice that Propositions
\ref{prop1} and \ref{prop2} below allow us to
construct an (exact) solution $u(y,
\lambda(\varepsilon), \varepsilon)$ of the
Rayleigh equation  hence giving an unstable mode
$${\tilde u}(x, z, t)= e^{ikz}u(\frac{x}{L_0},
kL_0,
\lambda(\varepsilon))e^{\frac{\sqrt{gk}}{\sqrt{\lambda(\varepsilon)}}t}$$ solution of the linearized Euler
equations.
Moreover, from Proposition (\ref{sc}), one has the following:
There exists $k\geq 1$, $\lambda(\varepsilon)$ and $u(y)$ such that $\varepsilon=kL_0$, $u$ solution of (\ref{eqrayleigh}), $\gamma(k, \varepsilon)=\sqrt{\frac{gk}{\lambda(\varepsilon)}}$, $\frac{\Lambda}{2}<\gamma(k, \varepsilon)<\Lambda$, $||u(y)||_{L^2}=1$, $u(0)>0$.\\
From the construction of this particular solution, we deduce a nonlinear result. For simplicity, in what follows, we will denote by $\gamma(k)$ the eigenvalue $\gamma(k, L_0)$.\\
From $u$, one deduces a solution $$U=\Re [(u_1, v_1, Q_1,
T_1)e^{ikz+\gamma(k)t}]=\Re [(u(x), -\frac{1}{ik}u'(x),
-\frac{\gamma(k)}{k^2}u'(x),
\frac{k_0(x)}{\gamma(k)}u(x))e^{ikz+\gamma(k)t}]$$ of the linearized
system. We thus consider a function $V^N=(0, 0, \frac{p_0}{\rho_0},
1)(L_0x)+\sum_{j=1}^N\delta^jV_j(x, y, t)$ satisfying
$(Emod)(V^N)=\delta^{N+1}R^{N+1}$, $V^N(x, z, 0)-(0, 0,
\frac{p_0}{\rho_0}, 1)(L_0x)=\delta U(x, z, 0)$. We also construct the
solution $V(x, y, t)$ of the Euler system such that $Emod(V)=0$ and
$V(x, z, 0)=(0, 0, \frac{p_0}{\rho_0}, 1)(L_0x)+\delta U(x, y,
0)$. Introduce finally $V^d(x, y, t)=V(x, y, t)-V^N(x, y, t)$. This
procedure constructs a solution of the nonlinear system.\\
We have the
\begin{theo}
\label{resNL}
\begin{enumerate}
\item There exists two constants $A$ and $C_0$, depending only on the properties of the Euler system, on the stationary solution and on the solution ${\hat u}(x)$, such that, for all $\theta<1$, for all $t\in ]0, \frac{1}{\gamma(k)}\ln\frac{\theta}{\delta C_0A}[$, one has the control of the approximate solution $V^N$ in $H^s$, namely
$$||T^N-1||_{H^s}+||{\vec u}^N||_{H^s}+||Q^N-q_0||_{H^s}\leq C\frac{\delta AC_0e^{\gamma(k)t}}{1-\delta AC_0e^{\gamma(k)t}}$$
and the leading order term of the approximate solution is the solution of the linear system as follows
$$||T^N-1||_{L^2}\geq \delta ||T_1(0)||_{L^2}e^{\gamma(k)t}- AC_0^2C_3\frac{e^{\gamma(k)t}}{1-\delta AC_0e^{\gamma(k)t}}$$
\item  There exists $N_0$ such that for any $N\geq N_0$, the function $V^d$ is well defined for  $t<\frac{1}{\gamma(k)}\ln \frac{1}{\delta}$ and satisfies the inequality
$$||V^d||\leq \delta^{N+1} e^{(N+1)\gamma(k)t}, \forall t\in [0, \frac{1}{\gamma(k)}\ln \frac{1}{\delta}[.$$
\item We have the inequality, for $\epsilon_0<C_0A$
$$||{\vec u}(\frac{1}{\gamma(k)}\ln\frac{\varepsilon_0}{C_0A\delta})||_{L^2}\geq\frac{\varepsilon_0}{2}||{\vec u}_1(0)||_{L^2}.$$
\end{enumerate}
\end{theo}
\vskip 1 cm

This paper is organized as follow. The sections 1, 2, 3 study the linear system and identify the behavior of the growth rate $\gamma(k)$ when $L_0\rightarrow 0$ by constructing the Evans function, and Section 4 constructs an approximate solution of the nonlinear system of Euler equations. \\
We identify in a first section
the family of solutions of (\ref{eqrayleigh}) which are bounded when
$y\rightarrow +\infty$ and we extend such solutions, for
$(\varepsilon, \lambda)$ in a compact ${\cal B}$, on $[\xi^{-1}(\varepsilon
R), +\infty[$, where $R$ is a constant depending only on ${\cal B}$
(Proposition \ref{prop1}). In the second section, for all $t_0>0$,
we calculate a solution of (\ref{eqrayleigh}) which is bounded on
$]-\infty,
-\frac{t_0}{\varepsilon}]$ (Proposition \ref{prop2}).\\
A solution $u$ of (\ref{eqrayleigh}) which is in $L^2(\mathbb R)$
goes to zero when $y\rightarrow +\infty$ as well as when
$y\rightarrow -\infty$. Moreover, as $\rho_0(x)$ is a $C^{\infty}$
function on $\mathbb R$, any solution $u$ of (\ref{eqrayleigh0}) is
also in $C^{\infty}$.\\ Notice that $\mbox{lim}_{\varepsilon \rightarrow 0}(-\varepsilon \xi^{-1}((\varepsilon R)^{\frac{1}{\nu}}))=\frac{1}{\nu R}$, from which 
one deduces that there exists $t_0$ such that $0<t_0<\frac12\mbox{lim}_{\varepsilon \rightarrow 0}(-\varepsilon \xi^{-1}((\varepsilon R)^{\frac{1}{\nu}}))$.\\ The regions $]-\infty, -\frac{t_0}{\varepsilon}]$ and $[\xi^{-1}((\varepsilon R)^{\frac{1}{\nu}}), +\infty[$ overlap and $$[\xi^{-1}((\varepsilon R)^{\frac{1}{\nu}}), -\frac{t_0}{\varepsilon}]\subset [-\frac{3}{4\varepsilon \nu R}, -\frac{1}{2\varepsilon \nu R}].$$
Hence the solution $u$ belongs to the family of
solutions described in proposition \ref{prop1}
(of the form $C_*u_+(y, \varepsilon)$) and belongs to the
family of solutions described in proposition
\ref{prop2} (of the form
$C_{**}U(-\varepsilon y, \varepsilon)$), that is

$$\bs\ba{l}u(y)=C_*u_+(y, \varepsilon), y\geq \xi^{-1}((\varepsilon R)^{\frac{1}{\nu}})\cr
u(y)=C_{**}U(-\varepsilon y, \varepsilon), y<-\frac{t_0}{\varepsilon}\ea\es$$
From the continuity of $u$ and of $u'$, one deduces that, for all $y_\perp\in [-\frac{3}{4\varepsilon \nu R}, -\frac{1}{2\varepsilon \nu R}]$ (corresponding to $t_\perp=-\varepsilon y_\perp\in [\frac{1}{2 \nu R}, \frac{3}{4 \nu R}]$), we have $C_*u_+(y_\perp, \varepsilon)=C_{**}U(t_\perp, \varepsilon)$, $C_*\frac{d}{dy}u_+(y_\perp, \varepsilon)= -C_{**}\varepsilon U'(t_\perp, \varepsilon)$.\\
Introduce the
Wronskian (where $\varepsilon^{-1}$ has been added for normalization purposes)
$$\mathcal{W}(y)=\varepsilon^{-1}(u_+(y,
\varepsilon)\frac{d}{dy}(U(-\varepsilon y, \varepsilon))-\frac{d}{dy}(u_+(y, \varepsilon))U(-\varepsilon y, \varepsilon)).$$ It is zero at $y_\perp=-\varepsilon
t_\perp$. Conversely, if $\lambda$ and $\varepsilon$ are chosen such that the Wronskian is zero (in particular at a point $y_\perp=-\frac{t_\perp}{\varepsilon}$), the function
\be\label{solution}{\tilde u}(y)=\bs\ba{l}C_{**}U(-\varepsilon y, \varepsilon), y\leq y_\perp\cr
C_{**}\frac{U(-\varepsilon y_\perp, \varepsilon)}{u_+(y_\perp,
\varepsilon)}u_+(y, \varepsilon), y\geq y_\perp\ea\es\ee is, thanks
to the Cauchy-Lipschitz theorem, a solution of (\ref{eqrayleigh0}).
Moreover, it belongs to $L^2(\mathbb R)$ owing to the properties of
$u_+$ and of $U$.\\In Section \ref{sec40}, we compute the function
$\mathcal{W}$. As $U$ and $u_+$ are solutions of the Rayleigh
equation, which rewrites
$$\frac{d^2}{dy^2}(u_+(y, \varepsilon))=-\frac{\xi'(y)}{\xi(y)}\frac{du_+}{dy}+(\varepsilon^2-\varepsilon \lambda \frac{\xi'(y)}{\xi(y)})u_+(y, \varepsilon)$$
the function $\mathcal{W}$ is solution of $\frac{d}{dy}\mathcal
W=-\frac{\xi'(y)}{\xi(y)}\mathcal W$, which implies the equality
\be\label{wr}\xi(y)\mathcal{W}(y)=\xi(y_0)\mathcal{W}(y_0) \mbox{
for all }y, y_0\ee This Wronskian can be computed for $y_\perp\in
[-\frac{3}{4\varepsilon \nu R}, -\frac{1}{2\varepsilon \nu R}]$
using the expressions obtained for $U$ and $u_+$. We prove that it
admits a unique root for $0<\varepsilon<\varepsilon_0$ and $\lambda$
in a fixed compact, and we identify the expansion of this root in
$\varepsilon$, hence proving Theorem \ref{HG0}. Precise estimates of this solution are given in Section 3. \\
In Section 4, after proving a $H^s$ result on a general solution of the linear system (taking into account a mixing of modes), we calculate all the terms $V_j$ of the expansion of the approximate solution, the perturbation of order $\delta$ being an eigenmode with a growth rate $\gamma\in ]\frac{\Lambda}{2}, \Lambda[$, where $\Lambda^2=\mbox{max}k_0(x)\frac{g}{L_0}$.
 \section{Construction of the family of bounded solutions in the dense region.} \label{sec2}  The
system (\ref{infinity}) writes $\frac{d}{dy}{\vec U}_+=\varepsilon
M_0(\xi(y), \lambda){\vec U}_+$. When $y\rightarrow
+\infty$, the matrix converges exponentially
towards $M_0(1, \lambda)$, which eigenvalues are
$0$ and $2$, of associated eigenvectors
$(1,\lambda-1)$ and $(1, \lambda+1)$.\\ 
It is classical that

\begin{lemma}
There
exists a unique solution $(U_+, V_+)$ of
(\ref{infinity}) which limit at $y\rightarrow
+\infty$ is $(1, \lambda -1)$.
Moreover, there exists $\xi_0>0$ such that this
solution\footnote{It can also be shown that there
exists a unique solution $({\tilde U}, {\tilde
V})$ such that $({\tilde U}, {\tilde
V})e^{-2\varepsilon y}\rightarrow (1, \lambda
+1)$} admits an analytic expansion in
$\varepsilon$ for $\xi(y)\in [\xi_0, 1[$.
\end{lemma}
The proof of this result is for example a consequence of Levinson
\cite{Lev}.\\
 The aim of this section is to
express precisely the coefficients of this
expansion when $\xi(y)\rightarrow 0$ and to deduce that one can extend the expression obtained for $\xi \in [\xi(\varepsilon R), \xi_0]$.\\
We consider, in what follows, the change of variable

\be \label{zeta}\zeta = \frac{\varepsilon}{\xi(y)^{\nu}}.\ee
 We prove in
this section the
\begin{prop} Let $K$ be a
compact set and $\lambda\in K$. There exists
$\varepsilon_0>0$ and $R>0$ such that, for
$0<\varepsilon<\varepsilon_0$, the family of
solutions of (\ref{infinity}) which is bounded
when $y\rightarrow +\infty$ is
characterized\footnote{a general solution is
$K_+(U, V)$ where $K_+$ is a constant}, for $y$
such that $\xi(y)\geq(\varepsilon
R)^{\frac{1}{\nu}}$, by
$$\bs\ba{l}U_+(y, \varepsilon)=1+\frac{(1-\xi(y))(1-\lambda)}{\xi(y)}\zeta A(\zeta, \varepsilon)\cr
V_+(y, \varepsilon)=\lambda -1+(1-\lambda)(1-\xi)\zeta B(\zeta,
\varepsilon).\ea\es$$ \label{prop1}
The associated solution of (\ref{eqrayleigh}) is $u_+(y, \varepsilon)=U_+(y, \varepsilon)e^{-\varepsilon y}$.
\end{prop}
\paragraph{Proof of Proposition \ref{prop1}}
We write the analytic expansion in $\varepsilon$:
$$U=1+\sum_{j\geq 1} \varepsilon^ju_j,V=\lambda -1+\sum_{j\geq 1} \varepsilon^jv_j.$$
We deduce, in particular,
$$\bs\ba{l}\frac{du_1}{dy}=\frac{\lambda -1}{\xi(y)}(1-\xi(y))\cr
\frac{dv_1}{dy}=(\lambda^2-1)(1-\xi(y))\ea\es$$
hence assuming $u_1, v_1\rightarrow 0$ when
$\xi\rightarrow 1$ (which is equivalent to
dividing the solution by its limit when
$\xi\rightarrow 1$) we get
$$\bs\ba{l}u_1=\frac{1-\lambda}{\nu+1}\frac{1-\xi^{\nu+1}}{\xi^{\nu+1}}\cr
v_1=\frac{1-\lambda^2}{\nu}\frac{1-\xi^{\nu}}{\xi^{\nu}}.\ea\es$$
The following recurrence system for $j\geq 1$
holds:
\be\label{coeffrecu}\bs\ba{l}\frac{du_{j+1}}{dy}=\frac{1}{\xi}(v_j-(\lambda-1)\xi
u_j)\cr \frac{dv_{j+1}}{dy}=(\lambda
+1)(v_j-(\lambda -1)\xi u_j).\ea\es\ee Usual
methods for asymptotic expansions lead to the
estimates (which are not sufficient for the proof of Proposition \ref{prop1})
$$|u_j(y)|+ |v_j(y)|\leq \frac{MA^j}{\xi_0^{(\nu+1)j}}.$$
However, using the relation $1-\xi = \frac{{\dot
\xi}}{\xi^{\nu+1}}$, we obtain the following
estimates:
\begin{lemma}
\label{inequality} Let $\xi_0>0$ given. For all
$j\geq 1$, introduce $a_j$ and $b_j$, such that
$$u_j(y)=\frac{(1-\xi(y))(1-\lambda)}{\xi^{\nu j +1}}a_j(\xi(y)), v_j(y)=\frac{(1-\xi(y))(1-\lambda)}{\xi^{\nu j}}b_j(\xi(y)).$$
The functions $a_j$ and $b_j$ are bounded,
analytic functions of $\xi$, for $\xi\in [\xi_0,
1]$. They satisfy \be\label{esti}|a_j(\xi)|\leq
AR^j, |b_j(\xi)|\leq AR^j,\ee where $R$ depends
only on $\lambda$.
\end{lemma} We prove Lemma \ref{inequality} by recurrence. Assume that this relation is
true for $j$. We have the relations
$$\bs\ba{l} \frac{du_{j+1}}{dy}=(1-\lambda)(b_j-(\lambda-1)a_j)\frac{{\dot
\xi}}{\xi^{\nu(j+1)+2}}\cr
\frac{dv_{j+1}}{dy}= (1-\lambda)(\lambda+1)(b_j-(\lambda-1)a_j)\frac{{\dot
\xi}}{\xi^{\nu(j+1)+1}}\ea\es$$
from which we deduce, using the limit 0 at $\xi\rightarrow 1$
$$u_{j+1}(y)=(1-\lambda)\int_1^{\xi(y)}\frac{b_j(\eta)-(\lambda-1)a_j(\eta)}{\eta^{\nu(j+1)+2}}d\eta$$
and
$$v_{j+1}(y)=(1-\lambda)(\lambda+1)\int_1^{\xi(y)}\frac{b_j(\eta)-(\lambda-1)a_j(\eta)}{\eta^{\nu(j+1)+1}}d\eta$$
We thus deduce that $\xi^{\nu(j+1)}v_{j+1}(y)$ and $\xi^{\nu(j+1)+1}u_{j+1}(y)$ are bounded functions when $\xi\in ]0, 1]$. Moreover, if we assume $|b_j|\leq AR^j$ and $|a_j|\leq AR^j$, then
$$\ba{l}|u_{j+1}|\leq AR^j|1-\lambda|(|\lambda -1|+1)\int_{\xi}^1\frac{d\eta}{\eta^{\nu(j+1)+2}}\cr |v_{j+1}|\leq AR^j|1-\lambda||\lambda +1|(|\lambda-1|+1)\int_{\xi}^1\frac{d\eta}{\eta^{\nu(j+1)+1}}.\ea$$
We end up with
$$\ba{l}|u_{j+1}|\leq |\lambda -1|AR^j\frac{(|\lambda -1|+1)}{\xi^{\nu(j+1)+1}}\frac{1-\xi^{\nu(j+1)+1}}{\nu(j+1)+1},\cr
 |v_{j+1}|\leq |\lambda -1|AR^j|\lambda +1|\frac{(|\lambda-1|+1)}{\xi^{\nu(j+1)}}\frac{1-\xi^{\nu(j+1)}}{\nu(j+1)}.\ea$$
As $\frac{1-\xi^a}{a}\leq 1-\xi, \xi\in [0, 1]$,
we get $|u_{j+1}|\leq |\lambda -1|AR^j\frac{(|\lambda
-1|+1)(1-\xi(y))}{\xi^{\nu(j+1)+1}}$,
$|v_{j+1}|\leq AR^j|\lambda -1||\lambda
+1|\frac{(|\lambda-1|+1)(1-\xi(y))}{\xi^{\nu(j+1)}}$.
Consider \be R_\lambda=(|\lambda
-1|+1)\mbox{max}(1, |\lambda +1|). \ee The
previous inequalities become $$|u_{j+1}|\leq
AR_\lambda^{j+1}\frac{(1-\xi(y))|\lambda -1|}{\xi^{\nu(j+1)+1}},
|v_{j+1}|\leq AR_\lambda
^{j+1}\frac{(1-\xi(y))|\lambda -1|}{\xi^{\nu(j+1)}},$$ hence
we proved the inequality for $j+1$.\\ The
inequality is true for $j=1$, hence the end of
the proof of Lemma \ref{inequality}, where we may choose the
value of $R$ for $\lambda\in [\frac12, \frac32]$
as $R=\frac{15}{4}$. Finally we have the
equalities, for all $y$ such that $\xi(y)\geq
\xi_0$:
$$\bs\ba{ll}U_+(y, \varepsilon)&=1+\frac{(1-\xi(y))(1-\lambda)}{\xi(y)}\sum_{j\geq 1}a_j(\xi(y))(\frac{\varepsilon}{(\xi(y))^{\nu}})^j\cr
&=
1+\frac{(1-\xi(y))(1-\lambda)}{\xi(y)}(\frac{\varepsilon}{(\xi(y))^{\nu}})\sum_{j\geq
0}a_{j+1}(\xi(y))(\frac{\varepsilon}{(\xi(y))^{\nu}})^j\cr
V_+(y, \varepsilon)&=\lambda -1+(1-\lambda)(1-\xi(y))\sum_{j\geq
1}b_j(\xi(y))(\frac{\varepsilon}{(\xi(y))^{\nu}})^j\cr
&=\lambda
-1+(1-\lambda)(1-\xi(y))(\frac{\varepsilon}{(\xi(y))^{\nu}})\sum_{j\geq
0}b_{j+1}(\xi(y))(\frac{\varepsilon}{(\xi(y))^{\nu}})^j.\ea\es$$
Using the
estimates (\ref{esti}) and the change of variable (\ref{zeta}), for $\zeta<R^{-1}$ the
series $\sum
a_j(\frac{\varepsilon^{\frac{1}{\nu}}}{\zeta^{\frac{1}{\nu}}})\zeta^j$
is normally convergent and the following
functions are well defined
$$\bs\ba{l}{\tilde U}(y, \varepsilon)=1+\frac{(1-\lambda)(1-\xi(y))}{\xi(y)}\zeta A(\zeta, \varepsilon)\cr
{\tilde V}(y, \varepsilon)=\lambda -1+(1-\lambda)(1-\xi(y))\zeta
B(\zeta, \varepsilon).\ea\es$$ It is
straightforward to check that ${\tilde U}$ and
${\tilde V}$ solve system (\ref{infinity}) and
that we have, for $\xi(y)\geq \xi_0$,
$\zeta(\xi)\leq \frac{\varepsilon}{\xi_0^{\nu}}$,
hence for\footnote{Note that these inequalities
depend on a given arbitrary $\xi_0>0$.}
$\varepsilon
<\varepsilon_0=\frac{\xi_0^{\nu}}{2R}$ and
$\xi(y)\geq\xi_0$ we have ${\tilde U}(y,
\varepsilon)=U_+(y, \varepsilon)$ and ${\tilde
V}(y, \varepsilon)=V_+(y, \varepsilon)$. We
extended the solution constructed for $\xi(y)\in
[\xi_0, 1[$ to the region $\zeta<\frac{1}{R}$.
This proves Proposition \ref{prop1}.
\section{The solution in the low density region}
\label{sec3}
\subsection{Construction of the bounded solution}In this section, we obtain the
family of solutions of (\ref{eqrayleigh}) bounded
by $|y|^Ae^{\varepsilon y}$ when $y\rightarrow
-\infty$, that is in the low density region
$\xi\rightarrow 0$. Introduce the new variable $t=-\varepsilon y$. Commonly, I call this solution the hypergeometric solution, because it has been observed that, in the model case $\rho_0(x)=(-x-1)^{-\frac{1}{\nu}}$ studied in \cite{CCLARA} as well as in \cite{Gon}, the Rayleigh equation rewrites as the hypergeometric equation.\\
Introduce 
$$\tau(s, \varepsilon)=-\frac{d}{ds}(\xi(-\frac{s}{\varepsilon}))(\xi(-\frac{s}{\varepsilon}))^{-1}=\frac{\xi^{\nu}}{\varepsilon}(1-\xi)= \frac{1}{\nu s}+\varepsilon^{\frac{1}{\nu}}S'(s,\varepsilon).$$
We define the operators $R_\varepsilon$, $K_\varepsilon$ and ${\tilde K}^{\lambda}_\varepsilon$ through

\be\label{operateurR}R_\varepsilon(g)(s)=[\int_s^{\infty}\tau(y, \varepsilon)e^{-2y}(\xi(-\frac{y}{\varepsilon}))^{-\lambda}g(y, \varepsilon)dy]e^{2s}(\xi(-\frac{s}{\varepsilon}))^{\lambda},\ee
\be\label{operateurtildeK}K_\varepsilon(g)(t)=(1-\lambda){\tilde K}^{\lambda}_{\varepsilon}(g)(t)=\frac{1-\lambda^2}{4}\int_t^{+\infty}\tau(s, \varepsilon)R_\varepsilon(g)(s, \varepsilon)ds.\ee
These operators rewrite

$$R_\varepsilon(g)(s, \varepsilon)=\int_s^{+\infty}(\frac{1}{\nu y}+\varepsilon^{\frac{1}{\nu}}S'(y,\varepsilon))
e^{-2(y-s)}s^{-\frac{\lambda}{\nu}}y^{\frac{\lambda}{\nu}}\exp(\varepsilon^{\frac{1}{\nu}}\lambda
(S(y)-S(s)))g(y, \varepsilon)dy.$$
$$K_\varepsilon(g)(t,
\varepsilon)=\frac{1-\lambda^2}{4}\int_t^{+\infty}(\frac{1}{\nu
s}+\varepsilon^{\frac{1}{\nu}}S'(s,\varepsilon))R_\varepsilon(g)(s,
\varepsilon)ds.$$

We have the inequalities, for $g$ uniformly bounded, (and $\lambda <\nu$, which implies $\xi(-\frac{s}{\varepsilon})^{\nu-\lambda}\leq \xi(-\frac{t}{\varepsilon})^{\nu - \lambda}$ for $t\geq s$)
\be |R_\varepsilon(g)(s)|\leq ||g||_{\infty}[\int_s^{+\infty}\frac{1}{\varepsilon}\xi^{\nu-\lambda}(1-\xi)e^{-2y}dy]e^{2s}\xi^\lambda\leq ||g||_{\infty}\frac{\xi^{\nu}}{\varepsilon}\ee
\be\label{inegK}|K_\varepsilon(g)(t)|\leq \frac{|\lambda^2-1|}{4}||g||_{\infty}\int_t^{\infty}\tau(s, \varepsilon)\frac{\xi^{\nu}}{\varepsilon}ds\leq \frac{|\lambda^2-1|}{4\nu}||g||_{\infty}\frac{\xi^{\nu}}{\varepsilon}.\ee
Moreover, the following inequality is true:
\be\label{inegKp}|g(s, \varepsilon)|\leq C_p(\frac{\xi^{\nu}}{\varepsilon})^p\Rightarrow |K_\varepsilon(g)(t, \varepsilon)|\leq \frac{|\lambda^2-1|}{8\nu (p+1)}C_p(\frac{\xi^{\nu}}{\varepsilon})^{p+1}.\ee
In a similar way, we introduce
$$\ba{l}K_0^{\lambda}(g)(t)=\frac{1-\lambda^2}{4}\int_t^{+\infty}\frac{1}{\nu s}R_0^{\lambda}(g)(s)ds\cr
R_0^{\lambda}(g)(s)=\int_s^{+\infty}\frac{1}{\nu
y}e^{-2(y-s)}s^{-\frac{\lambda}{\nu}}y^{\frac{\lambda}{\nu}}g(y)dy.\ea$$
Let $\varepsilon_0>0$ be fixed and $0<\varepsilon<\varepsilon_0$. Under suitable assumptions on $g$ (we can for example consider $g$ in $C^{\infty}([t_0, +\infty[)$ such that $|\partial^pg|\leq C_py^{\alpha -p}$ for all $p$), the operators $K_\varepsilon$, $R_\varepsilon$, $K_0$, $R_0$ are well defined. Moreover, one proves that
\be\label{volterraepsilon}g(t, \lambda, \varepsilon)=\sum_{n\geq 0}K_\varepsilon^{(n)}(1)(t, \varepsilon)\ee
\be\label{volterra}g_0(t, \lambda)=\sum_{n\geq 0}K_0^{(n)}(1)(t)\ee
are normally converging series on $[t_0, +\infty[$, and that we have:
\be g=1+K_\varepsilon(g), g_0=1+K_0(g_0).\ee 
Moreover, we know that $g$ is defined on $\mathbb R$, because the series $\sum \frac{(|\lambda^2-1|A)^p}{p!}(\frac{\xi^{\nu}}{\varepsilon})^p$ converges and is majorated by $\exp(|\lambda^2-1|A\frac{\xi^{\nu}}{\varepsilon})$, from the inequality (\ref{inegKp}). We obtain the inequalities
\be \label{majog}|g_0(t, \lambda)|\leq \exp(\frac{|\lambda^2-1|}{4\nu^2 t}), |g(t, \lambda,\varepsilon)|\leq \exp(\frac{|\lambda^2-1|}{8\nu}\zeta^{-1}).\ee
We cannot thus consider the limit $\zeta\rightarrow 0$ in the equalities containing $g$ as (\ref{majog}).\\
We shall assume that $\lambda$ belongs to a compact set and that $\lambda\geq\frac12$.
We prove
\begin{prop} \label{prop2}Let $g$ be defined through (\ref{volterraepsilon}). The family
of solutions of the system (\ref{hyp}) on $(u, w)$ which is
bounded by $|y|^Ae^{\varepsilon y}$ when
$y\rightarrow -\infty$ is given by $$u(y, \varepsilon)=C(F(t, \lambda, \varepsilon)+G(t, \lambda, \varepsilon)),
\xi(y)w(y, \varepsilon)=v(y, \varepsilon)=C\xi(y)[(\lambda-1)F(t, \lambda, \varepsilon)+(\lambda+1)G(t, \lambda, \varepsilon)]$$
where $C$ is a constant, $t\in [t_0, +\infty[$, $t=-\varepsilon y$ and
$F$ and $G$ are given by equalities
(\ref{fonctiong}) and (\ref{egaliteF}) below.\\
We have the estimates, for $t\in [t_0, \varepsilon[$
$$|g(t, \lambda, \varepsilon)-g_0(t, \lambda)|\leq C_0\varepsilon^{\frac{1}{\nu}}|g_0(t, \lambda)|$$
$$|u(-\frac{t}{\varepsilon}, \varepsilon)-u_0(-\frac{t}{\varepsilon}, \varepsilon)|\leq C_0\varepsilon^{\frac{1}{\nu}}|u_0(-\frac{t}{\varepsilon}, \varepsilon)|$$
$$|v(-\frac{t}{\varepsilon}, \varepsilon)-v_0(-\frac{t}{\varepsilon}, \varepsilon)|\leq C_0\varepsilon^{\frac{1}{\nu}}|v_0(-\frac{t}{\varepsilon}, \varepsilon)|$$
\end{prop}
\paragraph{proof} The system (\ref{hyp}) rewrites on $F$ and
$G$ given by Proposition \ref{prop2}: \be\label{fuchstoy}\bs\ba{l}F'(t, \lambda,
\varepsilon)=F(t, \lambda,
\varepsilon)-\frac12(\frac{1}{\nu
t}+\varepsilon^{\frac{1}{\nu}}S'(t,
\varepsilon))[(\lambda -1)F(t, \lambda,
\varepsilon)+(\lambda +1)G(t, \lambda,
\varepsilon)]\cr G'(t, \lambda,
\varepsilon)=-G(t, \lambda,
\varepsilon)+\frac12(\frac{1}{\nu
t}+\varepsilon^{\frac{1}{\nu}}S'(t, \varepsilon))
[(\lambda -1)F(t, \lambda, \varepsilon)+(\lambda
+1)G(t, \lambda, \varepsilon)].\ea\es\ee A non
exponentially growing solution of the system
(\ref{fuchstoy}) is obtained through the
following procedure. 
We denote by $g(t, \lambda, \xi)$ the function \be\label{fonctiong}\ba{ll}g(t,\lambda,
\varepsilon)&=G(t, \lambda,
\varepsilon)e^t(\xi(-\frac{t}{\varepsilon}))^{\frac{\lambda
+1}{2}}(\frac{\varepsilon}{\nu})^{-\frac{\lambda
+1}{2\nu}}\cr
&=G(t, \lambda,
\varepsilon)e^tt^{-\frac{1+\lambda}{2\nu}}\exp(-\varepsilon^{\frac{1}{\nu}}\frac{1+\lambda}{2}S(t,
\varepsilon)).\ea \ee
We first get, from the fact
that $F$ is bounded when $t\rightarrow +\infty$,
that

\be\label{egaliteF}\ba{ll}F(t,\lambda,
\varepsilon)e^{-t}t^{\frac{\lambda-1}{2\nu}}e^{\varepsilon^{\frac{1}{\nu}}\frac{\lambda-1}{2}S(t,
\varepsilon)} &=F(t,\lambda,
\varepsilon)e^{-t}(\xi(-\frac{t}{\varepsilon}))^{\frac{1-\lambda}{2}}(\frac{\varepsilon}{\nu})^{\frac{\lambda
-1}{2\nu}}\cr &= \frac{\lambda
+1}{2}\int_t^{+\infty}(\frac{1}{\nu
s}+\varepsilon^{\frac{1}{\nu}}S'(s,\varepsilon))s^{\frac{\lambda-1}{2\nu}}
e^{\varepsilon^{\frac{1}{\nu}}\frac{\lambda-1}{2}S(s,
\varepsilon)}e^{-s}G(s, \lambda,
\varepsilon)ds\cr
&=-\frac{\lambda +1}{2}\int_t^{+\infty}\xi^{-1}\frac{d}{ds}(\xi)g(s, \lambda, \varepsilon)e^{-2s}\xi^{-\lambda}(\frac{\varepsilon}{\nu})^{\frac{\lambda}{\nu}}ds\cr
=&-\frac{\lambda+1}{2}\int_t^{+\infty}\xi^{-1}\frac{d}{ds}(\xi)\xi^{\frac{1-\lambda}{2}}(\frac{\varepsilon}{\nu})^{\frac{\lambda-1}{2\nu}}ds.\ea\ee We deduce from the system (\ref{fuchstoy}) the equality

$$\frac{d}{dt}(G(t, \lambda, \varepsilon)e^tt^{-\frac{1+\lambda}{2\nu}}
\exp(-\varepsilon^{\frac{1}{\nu}}\frac{1+\lambda}{2}S(t,
\varepsilon)))= \frac{\lambda-1}{2}(\frac{1}{\nu
t}+\varepsilon^{\frac{1}{\nu}}S')e^tt^{-\frac{1+\lambda}{2\nu}}
\exp(-\varepsilon^{\frac{1}{\nu}}\frac{1+\lambda}{2}S(t,
\varepsilon))F(t, \lambda, \varepsilon).$$

Under the assumptions $g$ bounded and satisfies the condition
\be\label{limg}\mbox{lim}_{t\rightarrow \infty}g(t, \lambda, \varepsilon)=1\ee
one gets the equality
\be\label{equationg}g(t,\lambda,
\varepsilon)-1=K_\varepsilon(g)(t,
\varepsilon).\ee
Using the usual Volterra method and inequalities (\ref{inegK}), (\ref{inegKp}) and (\ref{majog}), we deduce that the only solution of (\ref{equationg}) satisfying assumptions (\ref{limg}) is given through (\ref{volterraepsilon}).\\
One gets $G$ through (\ref{fonctiong}) then $F$ thanks to
\be\label{Fdeg}F(t, \lambda, \varepsilon)e^{-t}\xi^{\frac{1-\lambda}{2}}(\frac{\varepsilon}{\nu})^{\frac{\lambda -1}{2\nu}}=(\frac{\varepsilon}{\nu})^{\frac{\lambda}{\nu}}\frac{\lambda+1}{2}\int_t^{\infty}\tau(s, \varepsilon)e^{-2s}\xi^{-\lambda}g(s, \lambda, \varepsilon)ds.\ee

 The first part of Proposition \ref{prop2} is proven.\\
Denote by $(u_0, w_0)$ the leading order term in $\varepsilon$ of $(u,
w)$ when $t$  and $\lambda$ are fixed. Introduce $F_0(t, \lambda)$ and $G_0(t, \lambda)$ through the equalities
$$u_0(t, \lambda)=F_0(t, \lambda)+G_0(t, \lambda), w_0(t, \lambda)=(\lambda-1)F_0(t, \lambda)+(\lambda+1) G_0(t, \lambda).$$
The functions $(F_0(t, \lambda), G_0(t, \lambda))$ are solution of
$$\bs\ba{l}\frac{dF_0}{dt}(t, \lambda)=F_0(t, \lambda)-\frac{\lambda -1}{2\nu t}F_0(t, \lambda)-\frac{\lambda +1}{2\nu t}
G_0(t, \lambda)\cr \frac{dG_0}{dt}(t,
\lambda)=-G_0(t, \lambda)+\frac{\lambda -1}{2\nu
t}F_0(t, \lambda)+\frac{\lambda +1}{2\nu t}G_0(t,
\lambda).\ea\es$$ The second part of Proposition \ref{prop2} comes from the following estimates on the operators $R_\varepsilon$ and $K_\varepsilon$, valid for
$\varepsilon\leq \varepsilon_0$ and $t\geq
t_0>0$:
\be\label{fuchsestimations}|R_\varepsilon(f)-R_0^{\lambda}(f)|\leq
C_1\varepsilon^{\frac{1}{\nu}}|R_0^{\lambda}(f)|,
|K_\varepsilon(g)-K_0^{\lambda}(g)|\leq
C_2\varepsilon^{\frac{1}{\nu}}|K_0^{\lambda}(g)|,\ee from which we deduce the
uniform estimates for $g$ given by (\ref{fonctiong}) solution of
(\ref{equationg})
\be |g(t,\lambda, \varepsilon)-g_0(t,
\lambda)|\leq
C_3\varepsilon^{\frac{1}{\nu}}|g_0(t, \lambda)|,
t\geq t_0, \varepsilon\leq \varepsilon_0\ee
because the Volterra series associated with $K_0$
is normally convergent in $[t_0, +\infty[$. This ends the proof of Proposition \ref{prop2}.\\
{\it Note that the previous
estimates, as well as the behavior of the
solution and the operator $R_0$, are valid only
for $t_0>0$, because, for example,
$R_0^{\lambda}(1)(s)\simeq 1$ when $s\rightarrow
0$. The integral defining $K_0^{\lambda}$ is
nevertheless convergent at $+\infty$, because for
$t\geq t_0$ we have the equality
$$2\nu sR_0^{\lambda}(1)(s)=
1-\int_s^{\infty}\frac{1}{y}(\frac{y}{s})^{\frac{1}{\nu}}e^{-2(y-s)}dy.$$}
\subsection{Construction of the hypergeometric solution for $\varepsilon=0$} 
We prove in this Section
\begin{lemma}
The solution $(F_0(t, \lambda), G_0(t, \lambda))$ constructed through (\ref{fonctiong}), (\ref{equationg}), (\ref{Fdeg}) for $\varepsilon=0$ is given by
$$\bs\ba{l}F_0(t, \lambda)=e^{-t}(U_0(t, \lambda)+\frac12\frac{dU_0}{dt}(t, \lambda))\cr
G_0(t, \lambda)=e^{-t}(U_0(t, \lambda)-\frac12\frac{dU_0}{dt}(t, \lambda))\ea\es$$
where $U_0(t)=2^{-\frac{\lambda+1}{2\nu}}U(-\frac{1+\lambda}{2\nu}, -\frac{1}{\nu}, 2t)$
the function $U(a, b, T)$ being the Logarithmic Kummer's solution of the confluent hypergeometric equation (see \cite{AS}).
\end{lemma}
This allows to obtain the limit of the $(F_0(t, \lambda),
G_0(t,\lambda))$ for $t\rightarrow 0$.\\
 The equation satisfied by
$U_0(t, \lambda)=u_0(t)e^{t}$  is\be
tU_0''-(2t+\frac{1}{\nu})U_0'
+\frac{\lambda+1}{\nu}U_0=0.\label{kummer}\ee Introducing
$T=2t$, we recognize (see \cite{AS}) the equation
for hypergeometric confluent functions for
$b=-\frac{1}{\nu}$ and
$a=-\frac{1+\lambda}{2\nu}$:
$$T\frac{d^2U_0}{dT^2}-(\frac{1}{\nu} +T)\frac{dU_0}{dT}+\frac{1+\lambda}{2\nu}U_0=0.$$
The family of solutions of this Kummer's equation is generated by two functions $M(a, b, T)$ and $U(a, b, T)$. Note that $T^{1-b}M(1+a-b, 2-b, T)$ is also a solution of (\ref{kummer}), independant of $M(a, b, T)$, hence $U(a, b, T)$ can be expressed using $M(a, b, T)$ and $T^{1-b}M(1+a-b, 2-b, T)$. The family of solutions of (\ref{kummer}) which go to zero when $T\rightarrow +\infty$ is generated by $U(a, b, T)$, called the logarithmic solution. It is the subdominant solution of the hypergeometric equation.\\ The expression of the subdominant solution $U(a, b, T)$ is the following:
$$U(a, b, T)=\frac{\pi}{\sin \pi b}[\frac{M(a, b, T)}{\Gamma(1+a-b)\Gamma(b)}-T^{1-b}\frac{M(1+a-b, 2-b, T)}{\Gamma(a)\Gamma(2-b)}]$$
where $\Gamma$ is the usual Gamma function ($\Gamma(s)=\int_0^{\infty}t^se^{-t}dt$).
The relation between $U(a, b, 0)$ and $U'(a,b, 0)$ characterize the subdominant solution of the ordinary differential equation, and this particular solution has been chosen through the limit\footnote{$\frac{\Gamma(1-b)}{\Gamma(1+a-b)}=\frac{\pi}{\sin
\pi b\Gamma(b)\Gamma(1+a-b)}$} when $z\rightarrow +\infty$:
\be
\label{kummer0}U(a, b,
0)=\frac{\Gamma(1-b)}{\Gamma(1+a-b)},
\mbox{lim}_{z\rightarrow +\infty}z^aU(a, b,
z)=1.\ee
As we imposed that $g(t,
\varepsilon)\rightarrow 1$ when $t\rightarrow
+\infty$, we get that $G_0(t,
\lambda)e^tt^{-\frac{\lambda+1}{2\nu}}\rightarrow
1$ when $t\rightarrow +\infty$ and that there exists a
constant ${\tilde C}$ such that $F_0(t,
\lambda)e^tt^{1-\frac{\lambda+1}{2\nu}}\rightarrow
{\tilde C}$ when $t\rightarrow +\infty$. Hence
$$(F_0(t, \lambda)+G_0(t,
\lambda))e^tt^{-\frac{\lambda+1}{2\nu}}\rightarrow
1.$$
As $T^aU(a, b, T)\rightarrow 1$, we get that $t^{-\frac{1+\lambda}{2\nu}}U(-\frac{1+\lambda}{2\nu}, -\frac{1}{\nu}, 2t)\rightarrow 2^{\frac{\lambda +1}{2\nu}}$. We thus obtain the equality
$$t^{-\frac{1+\lambda}{2\nu}}U_0(t, \lambda)=t^{-\frac{1+\lambda}{2\nu}}e^t(F_0(t, \lambda)+G_0(t, \lambda))= 2^{-\frac{\lambda +1}{2\nu}}t^{-\frac{1+\lambda}{2\nu}}U(-\frac{1+\lambda}{2\nu}, -\frac{1}{\nu}, 2t),$$
hence \be \label{fonctionbasehyper}U_0(t, \lambda)=2^{-\frac{\lambda
+1}{2\nu}}U(-\frac{1+\lambda}{2\nu}, -\frac{1}{\nu}, 2t).\ee
Introduce \be\label{const}C_0(\lambda)=U(-\frac{1+\lambda}{2\nu},
-\frac{1}{\nu}, 0)=-\frac{\pi}{\sin
\frac{\pi}{\nu}\Gamma(-\frac{1}{\nu})\Gamma(1+\frac{\lambda
-1}{2\nu})}= \frac{\Gamma(1+\frac{1}{\nu})}{\Gamma(1+\frac{\lambda
-1}{2\nu})}.\ee We get that $u_0(t)=2^{-\frac{\lambda
+1}{2\nu}}U(-\frac{1+\lambda}{2\nu}, -\frac{1}{\nu}, 2t)e^{-t}$. As
$w_0=\lambda u_0-\frac{du_0}{dt}= ((\lambda
+1)U_0-\frac{dU_0}{dt})e^{-t}$ one deduces \be\label{HCS}G_0(t,
\lambda)=(U_0(t, \lambda) - \frac{1}{2}\frac{dU_0}{dt}(t, \lambda))e^{-t}, F_0(t,
\lambda)=(U_0(t, \lambda)+ \frac{1}{2}\frac{dU_0}{dt}(t, \lambda))e^{-t}.\ee Using
\cite{AS} and (\ref{kummer}), we finally obtain \be\label{lim}G_0(t,
\lambda)\rightarrow 2^{-\frac{\lambda +1}{2\nu}}C_0(\lambda), F_0(t,
\lambda)\rightarrow 2^{-\frac{\lambda +1}{2\nu}}C_0(\lambda)\mbox{
when }t\rightarrow 0.\ee
We deduce the equality
$$U_0(0, \lambda)=2^{1-\frac{\lambda+1}{2\nu}}\frac{\Gamma(1+\frac{1}{\nu})}{\Gamma(1+\frac{1-\lambda}{2\nu})},
\mbox{lim}_{t\rightarrow +\infty} t^{-\frac{1}{\nu}-\frac{\lambda
-1}{2\nu}}U_0(t, \lambda)=2^{\frac{1}{\nu}+\frac{\lambda -1}{2\nu}}.$$

Note that we can deduce the expressions of $F_0+G_0$ and of $G_0$. We thus check that
\be (F_0+G_0)(t, \lambda)e^t=C_0(M(-\frac{1+\lambda}{2}, -\frac{1}{\nu}, 2t)-C_*(2t)^{\frac{1}{\nu}+1}M(1+\frac{1-\lambda}{2\nu}, 2+\frac{1}{\nu}, 2t))\ee
\be \ba{ll}e^tG_0(t, \lambda)=&C_0((M-M')(-\frac{1+\lambda}{2}, -\frac{1}{\nu}, 2t)-C_*(2t)^{\frac{1}{\nu}+1}(M-M')(1+\frac{1-\lambda}{2\nu}, 2+\frac{1}{\nu}, 2t))\cr
&-C_0C_*2^{\frac{1}{\nu}}(1+\frac{1}{\nu})t^{\frac{1}{\nu}}M(1+\frac{1-\lambda}{2\nu}, 2+\frac{1}{\nu}, 2t)).\ea\ee
We note that $(M-M')(-\frac{1+\lambda}{2}, -\frac{1}{\nu}, 0)=\frac{1-\lambda}{2}$. We deduce that $e^{t}(F_0+G_0)(0, \lambda)=2^{1-\frac{1+\lambda}{2\nu}}C_0(\lambda)$ and $e^tG_0(0, \lambda)=2^{-\frac{1+\lambda}{2\nu}}C_0(\lambda)(1-\lambda)$, hence $(\lambda-1)e^t(F_0+G_0)(0, \lambda)+2e^tG_0(0, \lambda)=0$.
 In the next Section, we combine the results of
Section \ref{sec2} and of Section \ref{sec3}.
\section{Precise calculus of the Evans function.}
\label{sec40}
The Wronskian is related to a function independant of the variable $t$, called the Evans function, introduced below in (\ref{evan}) and denoted by $Ev(\lambda, \varepsilon)$. In the present Section, we shall identify the leading order term in $\varepsilon$ of the Evans function, and all the terms of the form $\varepsilon^{\frac{1}{\nu}}(\lambda-1)$ of the Evans function. We shall finish by the calculation of the term of the form $\varepsilon^{\frac{2}{\nu}}$. More precisely, we prove
\begin{lemma}
\label{egaliteswr}
The function 

\be\label{evan}Ev(\lambda, \varepsilon)=\xi(y_0)\mathcal{W}(y_0)\ee is
independant of $y_0$. It is analytic in $\lambda$ and in
$\varepsilon^{\frac{1}{\nu}}$, $\varepsilon$. Moreover, one has $Ev(1,
\varepsilon)=2(\frac{\varepsilon}{\nu})^{\frac{1}{\nu}}$ and
$\partial_\lambda Ev(1,
0)=2^{1-\frac{1}{\nu}}\Gamma(1+\frac{1}{\nu})$. This function is called the Evans function of the equation (\ref{eqrayleigh}).
\end{lemma}
Using the expressions of $\frac{d}{dy}(U(-\varepsilon y, \varepsilon))$ and $\frac{d}{dy}u_+$, we have
$$\ba{ll}\varepsilon \mathcal{W}(y, \varepsilon)&=u_+(y, \varepsilon) (-\varepsilon \lambda U(-\varepsilon y, \varepsilon) + \varepsilon W(-\varepsilon y, \varepsilon))\cr
&- U(-\varepsilon y, \varepsilon)(-\varepsilon \lambda u_+(y, \varepsilon) + \frac{\varepsilon}{\xi(y)}v_+(y, \varepsilon))\cr
&=\frac{\varepsilon}{\xi(y)}(\xi(y)u_+(y, \varepsilon)W(-\varepsilon y, \varepsilon)-U(-\varepsilon y, \varepsilon)v_+(y, \varepsilon)).\ea$$
Hence we have the following constant function to study, which depends only on $\lambda, \varepsilon$:
$$Ev(\lambda, \varepsilon)=\xi(y)\mathcal{W}(y)= [\xi(y)u_+(y, \varepsilon)V(-\varepsilon y, \varepsilon)-\xi(y)v_+(y, \varepsilon)U(-\varepsilon y, \varepsilon)].$$
 We shall use the equalities, valid for all $y_0$ (and $t_0$) such that both solutions are defined (which means $y_0\in [-\frac{3}{4\varepsilon \nu R}, -\frac{1}{2\varepsilon \nu R}]$)
$$Ev(\lambda, \varepsilon)=\xi(y_0)\mathcal{W}(y_0)=(\xi\mathcal{W})(-\frac{t_0}{\varepsilon}).$$
We begin with the
\begin{lemma} The Evans function has an analytic expansion in $\lambda$, which coefficients depend analytically on $\varepsilon$ and $\varepsilon^{\frac{1}{\nu}}$.
\end{lemma}
For the precise study of the different terms of $Ev(\lambda, \varepsilon)$, we introduce
$$\xi=\xi(-\frac{t}{\varepsilon}), \zeta=\frac{\varepsilon}{\xi^\nu}, \zeta_0=\nu t= \zeta(t, 0), \mbox{ for }t\geq t_0>0.$$
We check that the function $Ev(\lambda, \varepsilon)$ is analytic in $\lambda$ and has an analytic expansion in $\varepsilon^{\frac{1}{\nu}}$ and $\varepsilon$ thanks to the equality
$$\sum_{p=0}^{[\nu]}\frac{1}{\xi^{\nu+1-p}}+\frac{1}{1-\xi}+\frac{1-\xi^{\nu - [\nu]}}{\xi^{\nu - [\nu]}(1-\xi)}= \frac{1}{\xi^{\nu+1}(1-\xi)}$$
which implies that the relation between $t$ and $\zeta$ is analytic in
$\varepsilon$ and $\varepsilon^{\frac{1}{\nu}}$.\\
Assume from now on $\lambda\geq \frac12$ and $\nu>2$ and replace $\xi(y)$ by $\varepsilon^{\frac{1}{\nu}}\zeta^{-\frac{1}{\nu}}$. 
Using this Lemma, there exists two functions $B_0(\varepsilon)$ and $C_0(\lambda, \varepsilon)$ such that
\be\label{defiB0C0}Ev(\lambda, \varepsilon)=Ev(1, \varepsilon)+B_0(\varepsilon)(\lambda -1)+C_0(\lambda, \varepsilon)(\lambda-1)^2.\ee
\paragraph{Direct relations}
Considering the limit in (\ref{lacamai}) for $\varepsilon=0$, we obtain
$$Ev(\lambda, 0)=(\lambda-1)e^t(F_0+G_0)(t, \lambda)[-1+(1-\lambda)\nu t A(\nu t, 0)+\nu tB(\nu t, 0)].$$
As this quantity is independant of $t$, we consider the limit when $t\rightarrow 0$, hence we deduce that \be \label{EV0}Ev(\lambda, 0)=-(\lambda -1)2^{1-\frac{\lambda +1}{2\nu}}C_0(\lambda).\ee Remark that this implies the {\bf identity}
\be (\lambda-1)e^t(F_0+G_0)(t, \lambda)[-1+(1-\lambda)\nu t A(\nu t, 0)+\nu tB(\nu t, 0)]=-2^{1-\frac{\lambda +1}{2\nu}}C_0(\lambda)(\lambda -1)\ee
which rewrites $$e^t(F_0+G_0)(t, \lambda)[-1+(1-\lambda)\nu t A(\nu t, 0)+\nu t B(\nu t, 0)]=-2^{1-\frac{\lambda +1}{2\nu}}C_0(\lambda).$$
In a similar way, we check  that, for $\lambda=1$, $U_+=1$ and $V_+=0$, and $g(t, 1, \varepsilon)=1$, which implies 
\be\label{G1}G(t, 1, \varepsilon)= e^{-t}(\frac{\nu }{\varepsilon})^{-\frac{1}{\nu}}\xi^{-1}=e^{-t}(\frac{\zeta}{\nu})^{\frac{1}{\nu}}\ee
from which one deduces
\be\label{EV1}Ev(1, \varepsilon)=2e^tG(t, 1, \varepsilon)\xi= 2(\frac{\varepsilon}{\nu})^{\frac{1}{\nu}}.\ee
From (\ref{defiB0C0}), the unique root $\lambda(\varepsilon)$ of $Ev(\lambda, \varepsilon)$ in the neighborhood of $\lambda =1$ satisfies
$$\lambda(\varepsilon)-1=-\frac{Ev(1, \varepsilon)}{B_0(\varepsilon)+C_0(\lambda(\varepsilon), \varepsilon)(\lambda(\varepsilon)-1)}.$$
The two first terms of the expansion of $\lambda(\varepsilon)-1$ in terms of $\varepsilon^{\frac{1}{\nu}}$ under the assumption $\nu>2$ are thus given through

$$\lambda(\varepsilon)-1=-\frac{Ev(1, \varepsilon)}{B_0(\varepsilon)+C_0(1,0)(\lambda(\varepsilon)-1)}+o(\varepsilon^{\frac{2}{\nu}}).$$
As $\lambda(\varepsilon)-1=-\frac{Ev(1, \varepsilon)}{B_0(0)}+o(\varepsilon^{\frac{1}{\nu}})$, we write
\be\label{devlambda}\ba{ll}\lambda(\varepsilon)-1&=-\frac{Ev(1, \varepsilon)}{B_0(\varepsilon)-C_0(1, 0)(B_0(0))^{-1}Ev(1, \varepsilon)}+o(\varepsilon^{\frac{2}{\nu}})\cr
&= -\frac{Ev(1, \varepsilon)}{B_0(\varepsilon)}-C_0(1, 0)(B_0(0))^{-3}(Ev(1, \varepsilon))^2+o(\varepsilon^{\frac{2}{\nu}}).\ea\ee

One is thus left with the calculus of $C_0(1, 0)$ and of $B_0(\varepsilon)$ up to the order 1.
For the computation of $B_0(\varepsilon)$, we need the behavior of the solutions of the overdense system for $\lambda=1$.

As in Section \ref{sec2}, we introduce $a_j(\xi)=a_j^0+\xi a_j^1+O(\xi^2)$ and $b_j(\xi)=b_j^0+\xi b_j^1+O(\xi^2)$. We recall that $\zeta A(\zeta, \varepsilon)=\sum_{j=1}^{\infty}a_j(\xi)\zeta^j$ and $\zeta B(\zeta, \varepsilon)=\sum_{j=1}^{\infty}b_j(\xi)\zeta^j$.
 Introduce 
$$u(\zeta)=\sum_{j\geq 1}\zeta^{j-1}b_j^0, v(\zeta)=\sum_{j\geq 1}\zeta^{j-1}jb_j^1, w(\zeta)=\sum_{j\geq 1}\zeta^{j-1}ja_j^0, k(\zeta)=\sum_{j\geq 1}\zeta^{j-1}ja_j^1.$$
\begin{lemma}
\label{relations} The following relations are true
$$\ba{l}e^t(F+G)(t, 1, \varepsilon)-e^t(F_0+G_0)(t, 1)=-\frac{\varepsilon^{\frac{1}{\nu}}}{\nu-1}+O(\varepsilon^{\frac{2}{\nu}})\cr
2e^tG(t, 1, \varepsilon)-2e^tG_0(t, 1)=-2\frac{\varepsilon^{\frac{1}{\nu}}}{\nu-1}+O(\varepsilon^{\frac{2}{\nu}})\cr
\zeta(t, \varepsilon)-\zeta_0(t)=-\xi \frac{\nu}{\nu-1}\zeta\cr
\zeta A(\zeta, \varepsilon)(1-\xi(y))-\zeta_0A(\zeta_0, 0)=\xi\zeta [k(\zeta)-w(\zeta)-\frac{\nu}{\nu-1}(\zeta w'(\zeta)+ w(\zeta))]+O(\varepsilon^{\frac{2}{\nu}})\cr
\zeta B(\zeta, \varepsilon)(1-\xi(y))-\zeta_0B(\zeta_0, 0)=\xi \zeta[v(\zeta)-u(\zeta)-\frac{\nu}{\nu-1}(\zeta u'(\zeta)+u(\zeta))]+O(\varepsilon^{\frac{2}{\nu}})\ea$$

For the computation of $C_0(1, 0)$, one has
$$C_0(1, 0)=-\mbox{lim}_{\lambda \rightarrow 1, t\rightarrow 0}\frac{e^t(F_0+G_0)(t, \lambda)-e^t(F_0+G_0)(t, 1)}{\lambda -1}.$$
\end{lemma}
From these two results, one obtains the following
\begin{prop}
\label{expan}
Introduce the function $R_0(t)=\frac{1}{2\nu}\ln t-{\tilde K}_0^1(1)(t)-\frac12 B_0(0)t^{-\frac{1}{\nu}}$, where ${\tilde K}^1_0(1)$ has been introduced in (\ref{operateurtildeK})and note that the terms $B_0$ and $C_0$ which have been introduced in (\ref{defiB0C0}) are calculated through
$$B_0(0)=-2\int_0^{+\infty}s^{\frac{1}{\nu}}e^{-2s}ds= -2^{-\frac{1}{\nu}}\Gamma(1+\frac{1}{\nu})$$
We have
$$B_0(\varepsilon)=B_0(0)+\frac{\varepsilon^{\frac{1}{\nu}}}{\nu -1}+2(\frac{\varepsilon}{\nu})^{\frac{1}{\nu}}\mbox{lim}_{t\rightarrow 0}R_0(t),$$
and
$$\ba{ll}C_0(1, 0)&=\int_0^{+\infty}s^{\frac{1}{\nu}}e^{-2s}ds-\frac{2}{\nu}\int_0^{+\infty}\ln s e^{-2s}ds+\frac{1}{\nu}\int_1^{+\infty}s^{\frac{1}{\nu}-1}e^{-2s}{\tilde K}_0^1(1)(s)ds\cr
&+\frac{B_0(0)}{2\nu}\int_0^1\frac{1-e^{-2s}}{s}ds+\frac{1}{\nu}\int_0^1s^{\frac{1}{\nu}-1}e^{-2s}[\frac{1}{2\nu}\ln s -R_0(s)]ds\ea$$
\end{prop}
Let us begin with the proof of Proposition \ref{expan}. We rewrite the Evans function as
\be\label{lacamai}\ba{ll}Ev(\lambda, \varepsilon)=&[(\lambda -1)e^t(F+G)(t, \lambda, \varepsilon)+2e^tG(t, \lambda, \varepsilon)](\xi(y)+(1-\lambda)(1-\xi)\zeta A(\zeta, \varepsilon))\cr
&- e^t(F+G)(t, \lambda, \varepsilon)(\lambda -1+(1-\lambda)(1-\xi(y))\zeta B(\zeta, \varepsilon)).\ea\ee
Remember that we have 
$$(\lambda-1)[B_0(\varepsilon)+C_0(\lambda, \varepsilon)(\lambda -1)]= Ev(\lambda, \varepsilon)-Ev(1, \varepsilon).$$
We thus deduce the equality
$$\ba{ll}B_0(\varepsilon)+C_0(\lambda, \varepsilon)(\lambda -1)&=\xi(y)\frac{2e^tG(t, \lambda, \varepsilon)-2e^tG(t, 1, \varepsilon)}{\lambda -1}\cr
&-(1-\xi(y))[e^t(F+G)(t, \lambda, \varepsilon)(1-\zeta B(\zeta, \varepsilon))+2e^tG(t, \lambda, \varepsilon)\zeta A(\zeta, \varepsilon)]\cr
&+(1-\lambda)(1-\xi(y))\zeta A(\zeta, \varepsilon)e^t(F+G)(t, \lambda, \varepsilon)\ea$$
Recall that $G(t, \lambda, \varepsilon)e^t=(\frac{\zeta}{\nu})^{\frac{\lambda +1}{2\nu}}g(t, \lambda, \varepsilon)$ and use $g(t, 1, \varepsilon)=1$.
We use also the relation (\ref{Fdeg}) to get 
\be\label{pourC}\ba{ll}(F+\frac{\lambda +1}{2\lambda}G)(t, \lambda, \varepsilon)e^{-t}(\frac{\zeta}{\nu})^{\frac{\lambda -1}{2\nu}}&=2\int_t^{\infty}\xi^{-\lambda}(\frac{\varepsilon}{\nu})^{\frac{\lambda}{\nu}}e^{-2s}g(s, \lambda, \varepsilon)ds\cr
&-\frac{\lambda +1}{2\lambda}\int_t^{\infty}\xi^{-\lambda}(\frac{\varepsilon}{\nu})^{\frac{\lambda}{\nu}}e^{-2s}[\frac{dg}{ds}-2(g-1)]ds\ea\ee
Note that we need two terms of $G$ and of $F+G$, and that we use $$\frac{dg}{ds}=(1-\lambda)\frac{d}{ds}({\tilde K}^{\lambda}_\varepsilon(g)),\quad g-1=(1-\lambda){\tilde K}^{\lambda}_\varepsilon(g).$$ This will contribute to the term in $C$. 
Rewrite the first term of (\ref{lacamai}) as
$$ \xi(y)\frac{2e^tG(t, \lambda, \varepsilon)-2e^tG(t, 1, \varepsilon)}{\lambda -1}=2(\frac{\varepsilon}{\nu})^{\frac{1}{\nu}}[\frac{(\frac{\zeta}{\nu})^{\frac{\lambda -1}{2\nu}}-1}{\lambda -1}-(\frac{\zeta}{\nu})^{\frac{\lambda -1}{2\nu}}{\tilde K}_\varepsilon^{\lambda}(g)].$$
Its limit when $\lambda$ goes to 1 is $2(\frac{\varepsilon}{\nu})^{\frac{1}{\nu}}[\frac{1}{2\nu}\ln(\frac{\zeta}{\nu})-{\tilde K}_\varepsilon^{1}(1)]$.
Hence we get the identity
\be\label{etoile}\ba{ll}B_0(\varepsilon)&=-(1-\varepsilon^{\frac{1}{\nu}}\zeta^{-\frac{1}{\nu}})[e^t(F+G)(t, 1, \varepsilon)(1-\zeta B_1(\zeta, \varepsilon))+2e^tG(t, 1, \varepsilon)\zeta A_1(\zeta, \varepsilon)]\cr
&+2(\frac{\varepsilon}{\nu})^{\frac{1}{\nu}}[\frac{1}{2\nu}\ln(\frac{\zeta}{\nu})-{\tilde K}_\varepsilon^{1}(1)]\ea\ee
and the right hand side is independant on $t$. Using Lemma \ref{relations}, we obtain
\be\label{eaza}\ba{ll}B_0(\varepsilon)&=2(\frac{\varepsilon}{\nu})^{\frac{1}{\nu}}[\frac{1}{2\nu}\ln(\frac{\zeta}{\nu})-{\tilde K}_\varepsilon^{1}(1)(t)]\cr
&+(1-\varepsilon^{\frac{1}{\nu}}\zeta^{-\frac{1}{\nu}})(-\frac{\varepsilon^{\frac{1}{\nu}}}{\nu-1})(1-\zeta B_1(\zeta, 0)+2\zeta A_1(\zeta, 0)+o(\varepsilon^{\frac{1}{\nu}}))\cr
&+(1-\varepsilon^{\frac{1}{\nu}}\zeta^{-\frac{1}{\nu}})[e^t(F_0+G_0)(t, 1)(1-\zeta B_1(\zeta, \varepsilon))+2e^tG_0(t, 1)\zeta A_1(\zeta, \varepsilon)]\ea\ee
from which one deduces
$$\ba{ll}B_0(\varepsilon)&=2(\frac{\varepsilon}{\nu})^{\frac{1}{\nu}}[\frac{1}{2\nu}\ln(\frac{\zeta}{\nu})-{\tilde K}_\varepsilon^{1}(1)]-\frac{\varepsilon^{\frac{1}{\nu}}}{\nu-1}(1-\zeta B_1(\zeta, 0)+2\zeta A_1(\zeta, 0)+o(\varepsilon^{\frac{1}{\nu}}))\cr
&-(1-\varepsilon^{\frac{1}{\nu}}\zeta^{-\frac{1}{\nu}})B_0(0)\cr
&-[e^t(F_0+G_0)(t, 1)(\zeta (B_1(\zeta, 0)- B_1(\zeta, \varepsilon))+2e^tG_0(t, 1)\zeta (A_1(\zeta, \varepsilon)-\zeta A_1(\zeta, 0))].\ea$$
Using the relations $G_0(t, 1)=(\frac{\zeta_0}{\nu})^{\frac{1}{\nu}}e^{-t}$ and $F_0(t, 1)=2e^t\int_t^{+\infty} s^{\frac{1}{\nu}}e^{-2s}ds$, one deduces that $G_0(t, 1)$and $F_0(t, 1)$ goes to a constant when $t\rightarrow 0$.\\
Hence one gets

\be\label{bzero}B_0(\varepsilon)=B_0(0)+(\frac{\varepsilon}{\nu})^{\frac{1}{\nu}}\mbox{lim}_{t\rightarrow 0}[\frac{1}{\nu}\ln t -2{\tilde K}_0^1(1)(t)-t^{-\frac{1}{\nu}}B_0(0))].\ee
The second part consists in the calculus of $C_0(1, 0)$.\\
Considering now $\varepsilon=0$ in (\ref{lacamai}), one obtains the two identities
$$B_0(0)=-e^t(F+G)(t, 1, 0)(1-\zeta_0B_1(\zeta_0, 0))-2e^tG(t, 1, 0)\zeta_0 A_1(\zeta_0, 0), \zeta_0=\nu t.$$
$$\ba{ll}B_0(0)+C_0(\lambda, 0)(\lambda -1)&=-e^t(F_0+G_0)(t, \lambda)(1-\zeta_0B(\zeta_0, 0))-2e^tG_0(t, \lambda)\zeta_0A(\zeta_0, 0))\cr
&-(\lambda -1)\zeta_0A(\zeta_0, 0)e^t(F_0+G_0)(t, \lambda).\ea$$
Hence
$$\ba{ll}C_0(\lambda, 0)(\lambda-1)=&-(\lambda -1)\zeta_0A(\zeta_0, 0)e^t(F_0+G_0)(t, \lambda)\cr
&+e^t(F_0+G_0)(t, 1)(\zeta B(\zeta, 0)-\zeta B_1(\zeta, 0))\cr
&+(1-\zeta_0 B(\zeta_0, 0))(e^t(F_0+G_0)(t, 1)-e^t(F_0+G_0)(t, \lambda))\cr
&+2e^tG_0(t, 1)(\zeta A_1(\zeta, 0)-\zeta A(\zeta, 0))\cr
&+\zeta_0 A(\zeta_0, 0)(2e^tG_0(t,1)-2e^tG_0(t,  \lambda)).\ea$$
We get (as we work for $\varepsilon=0$, we should write $\zeta_0$ but we drop this notation and we use $\zeta=\nu t$)
$$\ba{ll}C_0(\lambda, 0)=&-\zeta A(\zeta, 0)e^t(F_0+G_0)(t, \lambda)\cr
&-e^t(F_0+G_0)(t, 1)\zeta \frac{B_1(\zeta, 0)- B(\zeta, 0)}{\lambda -1}\cr
&-(1-\zeta B(\zeta, 0))\frac{e^t(F_0+G_0)(t, \lambda)-e^t(F_0+G_0)(t, 1)}{\lambda -1}\cr
&-2e^tG_0(t, \lambda)\zeta \frac{A(\zeta, 0)-\zeta A_1(\zeta, 0)}{\lambda -1}\cr
&-\zeta A(\zeta, 0)\frac{2e^tG_0(t, \lambda)-2e^tG_0(t, 1)}{\lambda -1}.\ea$$
In this equality, one only needs the value for $\lambda\rightarrow 1$, and it is independant of $\zeta$. We thus consider the limit when $\lambda\rightarrow 1$ and $\zeta\rightarrow 0$, hence one obtains
$$C_0(1, 0)=-\mbox{lim}_{\zeta\rightarrow 0, \lambda\rightarrow 1}\frac{e^t(F_0+G_0)(t, \lambda)-e^t(F_0+G_0)(t, 1)}{\lambda -1}.$$
Equality (\ref{pourC}) rewrites
$$\ba{c}(F+\frac{\lambda +1}{2\lambda}G)(t, \lambda, \varepsilon)e^{-t}(\frac{\zeta}{\nu})^{\frac{\lambda-1}{2\nu}}\cr
=\cr
 -\frac{\lambda +1}{2\lambda}\int_t^{\infty}(\varepsilon^{\frac{\lambda}{\nu}}\nu^{-\frac{\lambda}{\nu}}
\xi^{-\lambda})\frac{d}{ds}[e^{-2s}g(s, \lambda, \varepsilon)]ds\cr
+(\lambda -1)\frac{1+\lambda}{2\lambda}\int_t^{\infty}(\varepsilon^{\frac{\lambda}{\nu}}\nu^{-\frac{\lambda}{\nu}}\xi^{-\lambda})\frac{d}{ds}[{\tilde K}^{\lambda}_{\varepsilon}(g)e^{-2s}]ds\ea$$
Hence, considering the limit $\varepsilon\rightarrow 0$, one obtains
$$\ba{c}(F_0+\frac{\lambda +1}{2\lambda}G_0)(t, \lambda)e^{-t}(\frac{\zeta}{\nu})^{\frac{\lambda-1}{2\nu}}\cr
=\cr
 -\frac{\lambda +1}{2\lambda}\int_t^{\infty}(\frac{\zeta}{\nu})^{\frac{\lambda}{\nu}}\frac{d}{ds}[e^{-2s}]ds\cr
+(\lambda -1)\frac{1+\lambda}{2\lambda}\int_t^{\infty}(\frac{\zeta}{\nu})^{\frac{\lambda}{\nu}}\frac{d}{ds}[{\tilde K}^{\lambda}_{0}(g_0)e^{-2s}]ds\ea$$
The value for $\lambda=1$ is thus $(F_0+G_0)(t, 1)e^{-t}=2\int_t^{+\infty}(\frac{\zeta}{\nu})^{\frac{1}{\nu}}e^{-2s}ds$. Hence
$$\ba{c}e^{-t}[(\frac{1-\lambda}{2\lambda}G_0)(t, \lambda)(\frac{\zeta}{\nu})^{\frac{\lambda-1}{2\nu}}+(F_0+G_0)(t, \lambda)(\frac{\zeta}{\nu})^{\frac{\lambda-1}{2\nu}}-(F_0+G_0)(t, 1)]\cr
=\cr
 \frac{\lambda +1}{\lambda}\int_t^{\infty}(\frac{\zeta}{\nu})^{\frac{1}{\nu}}((\frac{\zeta}{\nu})^{\frac{\lambda-1}{\nu}}-1)e^{-2s}ds\cr
+(\lambda -1)\frac{1+\lambda}{2\lambda}\int_t^{\infty}\zeta^{\frac{\lambda}{\nu}}\frac{d}{ds}[{\tilde K}^{\lambda}_{0}(g_0)e^{-2s}]ds\ea$$
Dividing by $\lambda -1$, one deduces
$$\ba{c}e^{-t}[-\frac{1}{2\lambda}G_0(t, \lambda)(\frac{\zeta}{\nu})^{\frac{\lambda-1}{2\nu}}+\frac{(F_0+G_0)(t, \lambda)-(F_0+G_0)(t, 1)}{\lambda -1}(\frac{\zeta}{\nu})^{\frac{\lambda-1}{2\nu}}+(F_0+G_0)(t, 1)\frac{(\frac{\zeta}{\nu})^{\frac{\lambda-1}{2\nu}}-1}{\lambda -1}]\cr
=\cr
 \frac{\lambda +1}{\lambda}\int_t^{\infty}\frac{(\frac{\zeta}{\nu})^{\frac{\lambda}{\nu}}-1}{\lambda -1}e^{-2s}ds\cr
+\frac{1+\lambda}{2\lambda}\int_t^{\infty}(\frac{\zeta}{\nu})^{\frac{\lambda}{\nu}}\frac{d}{ds}[{\tilde K}^{\lambda}_{0}(g_0)e^{-2s}]ds\ea$$
We consider the limit when $\lambda\rightarrow 1$, and recalling that for $\varepsilon=0$ one has $\frac{\zeta}{\nu}=s$, denoting by $H(t, \lambda)=\frac{(F_0+G_0)(t, \lambda)-(F_0+G_0)(t, 1)}{\lambda -1}$, we obtain
$$\ba{c}e^{-t}[-\frac12 G_0(t, 1)+H(t, 1)+(F_0+G_0)(t, 1)\frac{1}{2\nu}\ln\frac{\zeta}{\nu}]\cr
=\cr
 2\int_t^{\infty}\ln s e^{-2s}ds\cr
+\int_t^{\infty}s^{\frac{1}{\nu}}\frac{d}{ds}[{\tilde K}^{1}_{0}(1)e^{-2s}]ds\ea$$
Using again the integration by parts on the last term hence one gets
$$\ba{c}e^{-t}[-\frac12 G_0(t, 1)+H(t, 1)+(F_0+G_0)(t, 1)\frac{1}{2\nu}\ln\frac{\zeta}{\nu}]\cr
=\cr
 2\int_t^{\infty}\ln s e^{-2s}ds
-t^{\frac{1}{\nu}}{\tilde K}^{1}_{0}(1)(t)e^{-2t}-\int_t^{\infty}\frac{d}{ds}(s^{\frac{1}{\nu}}){\tilde K}^{1}_{0}(1)e^{-2s}ds\ea$$
We notice that the function $R_0(t)=-\frac{1}{\nu}\ln t+2{\tilde K}_0^1(1)(t)+t^{-\frac{1}{\nu}}B_0(0)$ has a finite limit when $t$ goes to zero, according to (\ref{bzero}). We have the equality ${\tilde K}_0^1(t)=\frac12 R_0(t)+\frac{1}{2\nu}\ln t -\frac{1}{2}B_0(0)t^{-\frac{1}{\nu}}$. We deduce that
$$\ba{ll}\int_t^{1}\frac{d}{ds}(s^{\frac{1}{\nu}}){\tilde K}^{1}_{0}(1)e^{-2s}ds&= \int_t^{1}\frac{d}{ds}(\zeta^{\frac{1}{\nu}})[\frac12 R_0(s)+\frac{1}{2\nu}\ln s-\frac12 B_0(0)s^{-\frac{1}{\nu}}]e^{-2s}ds\cr
&=\nu^{\frac{1}{\nu}-1}\int_t^{1}s^{\frac{1}{\nu}-1}[\frac12 R_0(s)+\frac{1}{2\nu}\ln s-\frac12 B_0(0)s^{-\frac{1}{\nu}}]e^{-2s}ds\ea$$
In this last term, the only term which matters when $t\rightarrow 0$ is the term 
$$-\frac12 B_0(0)\nu^{\frac{1}{\nu}-1}\int_t^{1}s^{-1}e^{-2s}ds= -\frac12 B_0(0)\nu^{\frac{1}{\nu}-1}[\int_t^{1}\frac{e^{-2s}-1}{s}ds-\ln t].$$
Note that $B_0(0)=-2\int_0^{+\infty}s^{\frac{1}{\nu}}e^{-2s}ds$. One obtains
$$B_0(0)=-\int_0^{\infty}(\frac{a}{2})^{\frac{1}{\nu}}e^{-a}da=-2^{-\frac{1}{\nu}}\Gamma(1+\frac{1}{\nu}).$$

\subsection{Reduction of the Evans function}

\paragraph{Lower order terms}

Recall that the operator ${\tilde K}^1_\varepsilon(1)$ is defined through (\ref{operateurtildeK}). We prove the following lemma of reduction:

\paragraph{Proof of Lemma \ref{relations}}
 It is enough to prove that the relation giving $\zeta$ is
$$-\frac{t}{\varepsilon}=C-\frac{1}{\nu \xi^\nu}-\frac{1}{(\nu-1)\xi^{\nu-1}}-\xi^{2-\nu}R(\xi)$$
hence we deduce
$$t=-C\varepsilon+\frac{\zeta}{\nu}+\xi\frac{\zeta}{\nu-1}-R(\xi)\xi^2\zeta.$$
We thus obtain $t=\frac{\zeta_0}{\nu}$, hence
$$\frac{\zeta_0-\zeta}{\nu}=\xi\frac{\zeta}{\nu-1}+O(\xi^2)\zeta.$$
We deduce that $w(\zeta)-w(\zeta_0)=(\zeta - \zeta_0)w'(\zeta_0)+O((\zeta -\zeta_0)^2)$, hence $w(\zeta)-w(\zeta_0)-(\zeta - \zeta_0)w'(\zeta_0)=0(\varepsilon^{\frac{2}{\nu}})$ and $w(\zeta)-w(\zeta_0)-(\zeta - \zeta_0)w'(\zeta)=0(\varepsilon^{\frac{2}{\nu}})$.\\
We use $e^tG(t, \lambda, \varepsilon)=(\frac{\varepsilon}{\nu})^{\frac{\lambda+1}{2\nu}}\xi^{-\frac{\lambda+1}{2}}g(t, \lambda, \varepsilon)$, hence for $\lambda = 1$ we obtain
$$e^tG(t, 1, \varepsilon)=\eta^{-1}.$$
The equality giving $F(t, 1, \varepsilon)$ being
$$e^{-t}F(t, 1, \varepsilon)= \int_t^{+\infty}\tau(s, \varepsilon)e^{-2s}\eta(s, \varepsilon)^{-1}ds=-e^{-2t}\eta(t, \varepsilon)^{-1}+2\int_t^{+\infty}e^{-2s}\eta(s, \varepsilon)^{-1}ds,$$
one obtains
$$e^t(F+G)(t, 1,\varepsilon)-e^t(F_0+G_0)(t, 1)=2e^{2t}\int_t^{\infty}e^{-2s}(\frac{1}{\eta(s, \varepsilon)}-\frac{1}{\eta(s, 0)})ds.$$
Similarily
$$e^t(2G(t, 1, \varepsilon)-2G_0(t, 1))=\frac{2}{\eta(s, \varepsilon)}[1-\frac{\eta(s, \varepsilon)}{\eta(s, 0)}].$$
Using the relation $$\frac{1}{\eta(t, \varepsilon)^{\nu}}(1+\frac{\nu}{\nu-1}\varepsilon^{\frac{1}{\nu}}\eta(t, \varepsilon)+O(\varepsilon^{\alpha}))=\frac{1}{\eta(t, 0)^{\nu}}$$
one obtains 
$$\frac{\eta(t, \varepsilon)}{\eta(t, 0)}-1=\frac{1}{\nu-1}\varepsilon^{\frac{1}{\nu}}\eta(t, 0)+O(\varepsilon^{\alpha}).$$
This gives directly the two equalities of Lemma \ref{relations}.\\

\subsection{Limit for large $k$ of the growth rate}
Recall that was proven in \cite{HelLaf} the
following estimate on any value of $\gamma$ such
that there exists a solution of
(\ref{eqrayleigh}) associated with
$\lambda=\frac{gk}{\gamma^2}$ and
$\varepsilon=kL_0$:
$$\gamma\rightarrow
\Lambda=\sqrt{\frac{g}{L_0}}\sqrt{\frac{\nu^{\nu}}{(\nu+1)^{\nu+1}}}\mbox{
when }k\rightarrow +\infty.$$ If we compare with
(\ref{gamma}), one may see the difference between
the result for $L_0\rightarrow 0$ when $k$ is
fixed and the result for $L_0>0$ fixed and
$k\rightarrow +\infty$. Note for example that the
limit of
$$\frac{\sqrt{gk}}{\sqrt{1+(\frac{kL_0}{\nu})^{\frac{1}{\nu}}\Gamma(1+\frac{1}{\nu})}}$$
when $k\rightarrow +\infty$ is $+\infty$ because $\nu>1$. This is not surprising because
we did not get the lower order terms up to the order $\varepsilon$ of the expansion of $\lambda$.
Remark that the term in $\varepsilon$ comes from the terms in $\varepsilon$ in the functions
$A(\zeta, \varepsilon)$ and $B(\zeta, \varepsilon)$.\\
We have the following result (according to \cite{HelLaf})
\begin{lemma}
\label{existencenormalmode}
There exists $k_*>0$ such that, for all $k\geq k_*$, there exists a real $\gamma(k)$ and a non zero solution $u(x)e^{iky+\gamma(k)t}$ of the Rayleigh equation (\ref{eqrayleigh0}) such that 
$$\frac{\Lambda}{2}<\gamma(k)<\Lambda.$$
We have the following behavior of the eigenmode
$$||\rho_0^{\frac12}u||+||\rho_0^{\frac12}u'||+||u||+||u'||+||u''||<+\infty$$
\end{lemma}
As the result of this Lemma is important for the nonlinear analysis, we rewrite an idea of the proof, based on Remark 8.1 of \cite{HelLaf}. We denote by $L^2_{\rho_0^\frac12}$ the space of functions $u$ such that $\rho_0^{\frac12} u\in L^2(\RR)$.\\
{\bf Finding $\gamma$ is equivalent to finding 0 as an eigenvalue (in $L^2(\RR)$) of}
$$-\frac{1}{k^2}\rho_0^{-\frac12}\frac{d}{dx}(\rho_0\frac{d}{dx}\rho_0^{-\frac12})+1-\frac{g}{\gamma^2}k_0(x).$$
This operator rewrites $-\frac{1}{k^2}\frac{d^2}{dx^2}+1-\frac{g}{\gamma^2}k_0(x)+k^{-2}W_0(x)$
where $W_0(x)=\frac12 k'_0(x)+\frac14(k_0(x))^2$, which is bounded when$\frac{\rho''_0}{\rho_0}$ is bounded (or equivalently when $k'_0$ is bounded). We introduce the operator $Q=-\frac{1}{k^2}\rho_0^{-\frac12}\frac{d}{dx}(\rho_0\frac{d}{dx}\rho_0^{-\frac12})+1$, which is coercive, thanks to the Poincare estimates, for $k$ large enough. The eigenvalue problem rewrites
$$\frac{\gamma^2}{g}\in \sigma_p(Q^{-\frac12}k_0Q^{-\frac12}).$$
Under the (natural) hypothesis that $k_0$ has a nondegenerate minimum $L_0$, one deduces that for $k$ large enough one has at least a value of $\gamma(k)$ such that $L_0<\frac{g}{(\gamma(k))^2}<4L_0$ using usual results on semiclassical Schrodinger operators which potential has a well.\\ We thus constructed $v\in L^2(\RR)$ and $\gamma(k)$ such that $v$ is the eigenvector of $Q^{-\frac12}k_0Q^{-\frac12}$ associated with the eigenvalue $\frac{(\gamma(k))^2}{g}$.\\
To $v$ is associated a solution of (\ref{eqrayleigh0}) which is $u=\rho_0^{-\frac12}Q^{-\frac12}v$, $u'\in L^2_{\rho_0^{\frac12}}, u\in L^2_{\rho_0^{\frac12}}$. Remembering that $u$ solves
$$-u''+k^2u-k_0(x)u'-\frac{gk^2}{(\gamma(k))^2}u=0,$$
multiplying this equation by $u$ and integrating, one gets
$$\int (k^2u^2+(u')^2)dx= \int k_0(x)\rho_0^{-\frac12}u.[\frac{gk^2}{(\gamma(k))^2}\rho_0^{\frac12}u+\rho_0^{\frac12}u']dx$$
hence, using the hypothesis
$$\rho'_0\rho_0^{-\frac32}\leq M$$
one obtains (the norm on the Sobolev space $H^1$ is $||u||_1^2=\int (u')^2+k^2u^2dx$)
$$||u||_1\leq M[\frac{gk^2}{(\gamma(k))^2}||\rho_0^{\frac12}u||+||\rho_0^{\frac12}u'||]$$
hence a control on the $H^1$ norm of $u$ (instead of having the weight $\rho_0^{\frac12}$).\\
Moreover, as $u''=\frac{gk^2}{(\gamma(k))^2}k_0(x)u+k_0(x)u'-k^2u$, one deduces that $u''\in L^2$, and we have iteratively the control of $u$ in $H^s$ ($s\leq s_{max}$, according to the number of derivatives of $k_0$ that we consider).
\vfill
\eject
\section{Towards a non linear analysis}
We show in this Section that the result of Guo and Hwang \cite{GH} can be extended in our set-up, even if the density profile $\rho_0(x)$ does not satisfy the coercivity assumption (3) of \cite{GH}. The quantity $k_0(x)=\frac{\rho'_0(x)}{\rho_0(x)}$ plays a crucial role. It has a physical interpretation, being the inverse of a length: it is called the inverse of the density gradient scalelength. We need the assumptions
$$(H)\quad k_0(x) \mbox{ bounded }, k_0(x)\rho_0^{-\frac12} \mbox{ bounded}.$$ Note that $k_0$ bounded is fulfilled in the case studied by Guo and Hwang (where $\rho_0$ is bounded below), and in the case of the striation model (studied by R. Poncet \cite{poncet}) but is not automatically fulfilled by a profile such that $\rho_0(x)\rightarrow 0$ when $x\rightarrow -\infty$. However, for the particular case of the ablation front profile, we have $k_0(x)=L_0^{-1}\xi(\frac{x}{L_0})^{\nu}(1-\xi(\frac{x}{L_0}))$, hence it is bounded and belongs to $[0, L_0^{-1}\frac{\nu^{\nu}}{(\nu+1)^{\nu+1}}]$.\\

Before starting the proof of Theorem \ref{resNL}, which is rather technical, let us describe our procedure.\\
Firstly, we prove that the linear system reduces to an elliptic equation on the pressure, from which we obtain a general solution. We identify a normal mode solution of this system using the first part of the paper.\\
Once this normal mode solution $U$ is constructed, with suitable assumptions on the growth rate, one introduces a perturbation solution of the nonlinear system, which initial condition is $\delta  U\vert_{t=0}$ and an approximate solution $V^N$ of the non linear system which admits an expansion in $\delta^N$ up to the order $N$ with the same initial condition.\\
Using the Duhamel principle for the construction of the $j-$th term of the expansion in $\delta$ of $V^N$, one obtains a control of all the terms of $V^N$.\\
The natural energy inequalities are on the quantities $\rho_0^{\frac12}u^j$, $\rho_0^{\frac12}v^j$, $\rho_0^{-\frac12}p^j$, $\rho_0^{-\frac12}\rho^j$. We verify that the properties of $\rho_0(x)$ imply that we can {\bf deduce} inequalities on $u^j$, $v^j$, $\rho_0^{-1}p^j$ and $T^j$.\\
Note that we have, as a consequence of the method that we chose, a control in $t^se^{\Lambda t}$ of the $H^s$ norm of all solutions of the homogeneous linear system (with any initial condition $U(x, y, 0)$), and a control by $e^{j\gamma(k)t}$ (with no additional power in $t$) of the $H^s$ norm of the $j-$th term of the expansion.\\
\begin{rem}
When an initial value mixes eigenmodes, the $H^s$ norm of the solution  behaves as $t^se^{\Lambda t}$. If one starts from a pure eigenmode with $\frac{\Lambda}{2}<\gamma(k)<\Lambda$ the exponential behavior comes at most from the growth of the pure eigenmode.
\end{rem}

\subsection{Obtention of a solution of the linear system}
Consider the system
$$\bs\ba{l}\partial_t\sigma+\rho'_0v_1=f_0\cr
\rho_0\partial_tv_1+\partial_xp=\sigma g+f_1\cr
\rho_0\partial_tv_2+\partial_yp=f_2\cr
\partial_xv_1+\partial_yv_2=0\ea\es$$
We know that the relevant quantities are $\rho_0^{\frac12}v_{1, 2}$, $\rho_0^{-\frac12}\sigma$, and we denote these three quantities by $X, Y, \tau$.
To have the same behavior when $\rho_0\rightarrow 0$, consider $\psi$ such that, once $\psi$ is obtained, we revert to $v_1$ and $v_2$ using $v_1=-\partial_y(\rho_0^{-\frac12}\psi)$, $v_2=\partial_x(\rho_0^{-\frac12}\psi)$. Introduce
\be\label{inconnue auxiliaire}b=\rho_0^{-\frac12}[\partial_y(\rho_0v_1)-\partial_x(\rho_0v_2)].\ee
The system on $v_1, v_2, \sigma, p$ implies the two equations
\be\bs\ba{l}\partial_tb=g\partial_y\tau + \rho_0^{-\frac12}(\partial_yf_1-\partial_xf_2)\cr
\partial_t\tau+k_0(x)X=\rho_0^{-\frac12}f_0.\ea\es\ee
We obtain $\psi$ from $b$ through the elliptic equation
\be\Delta \psi -(\frac12k'_0+\frac14k_0^2)\psi=-b.\ee
We then revert to $X$ through the equality $X=-\partial_y\psi$.  Finally, the pressure $p$ is obtained through the elliptic equation$$\rho_0\partial_x(\rho_0^{-1}\partial_xp)+\partial^2_{y^2}p=\rho_0^{\frac12}[\rho_0^{\frac12}\partial_x(\rho_0^{-\frac12}\tau)g+\rho_0^{\frac12}\partial_x(\rho_0^{-1}f_1)+\rho_0^{-\frac12}f_2]$$
which rewrites 
\be \Delta p-k_0\partial_xp=\rho_0^{\frac12}[\partial_x\tau g-\frac12k_0\tau g)+\rho_0^{-\frac12}(\mbox{div}{\vec f}-k_0f_1)]\ee
Hence we solve the system
\be\bs\ba{l}\partial_t\tau=k_0\partial_y\psi(b)+\rho_0^{-\frac12}f_0\cr
\partial_tb=g\partial_y\tau + \rho_0^{-\frac12}(\partial_yf_1-\partial_xf_2)\cr
\tau(0)=\tau_0(x, y), b(0)=b_0(x, y)\cr
\Delta \psi -(\frac12k'_0+\frac14k_0^2)\psi=-b\cr
\Delta p-k_0\partial_xp=\rho_0^{\frac12}[\partial_x\tau g-\frac12k_0\tau g)+\rho_0^{-\frac12}(\mbox{div}{\vec f}-k_0f_1)]\ea\es\ee
which has the same properties as the system (13) of \cite{GH}, the
Poincare estimate being still valid. \\
From $b$ and $\tau$, one reverts to $X$ and $Y$, hence a solution of the system. Moreover, one checks that $(X, Y)\in L^2(\RR)$ (according to the energy equality), hence $X\in H^1(\RR)$ under the assumption $k_0$ bounded.

 \begin{prop}
 \label{controleQ}
 Under the hypotheseses (H), and under the hypothesis $h_j\in L^2, j=0, 1, 2$, the functions $u_1, v_1, T_1, p_1$ solution of
 $$\bs\ba{l}\partial_tT_1-k_0u_1=h_0\cr
 \rho_0\partial_tu_1+\partial_xp_1+\rho_0gT_1=h_1\cr
 \rho_0\partial_tv_1+\partial_yp_1=h_2\cr
 \partial_xu_1+\partial_yv_1=0\ea\es$$
 satisfies $u_1(t), v_1(t), T_1(t)\in L^2$ when it is true for $t=0$. Moreover, one has $\rho_0^{-1}p_1(t)\in L^2(\RR^2)$.
 \end{prop}
 \paragraph{Proof}
 The proof of this result follows two steps: first of all the assumption $k_0$ bounded implies that $\rho_0^{\frac12}u_1, \rho_0^{\frac12}v_1, \rho_0^{-\frac12}\nabla p_1, \rho_0^{\frac12}T_1$ belong to $L^2$. We thus multiply the equality $\partial_t{\vec u}_1+\rho_0^{-1}\nabla p_1+ T_1{\vec g}={\vec h}$ by $\nabla(\rho_0^{-1}p)$. We get, integrating in $x,y$:
 $$\int (\nabla q_1)^2+k_0(x)q_1\nabla q_1.{\vec e}_1+ T_1{\vec g}\nabla q_1=\int {\vec h}\nabla q_1$$
 from which one deduces
 $$||\nabla q_1||\leq \mbox{max}(k_0\rho_0^{-\frac12})||\rho_0^{\frac12}q_1||+g||T_1||_{\infty}+||{\vec h}||.$$
 It is then enough to use the Poincare estimate between $\rho_0^{\frac12}q_1$ and $\rho_0^{-\frac12}\nabla p_1$ to obtain the estimate on $\nabla q_1$, from which one deduces the estimate on $q_1$.\\
 Finally, from the estimate on $q_1$ and on $\nabla q_1$, multiplying the equation on the velocity by ${\vec u}_1$ and integrating, we get the Gronwall type inequality
 $$\frac{d}{dt}||{\vec u}_1||\leq C||q_1||_{H^1}+||{\vec h}||+g||T_1||_{\infty}$$
 hence a control on $||{\vec u}_1||$ on $[0, T]$ for all $t$ as soon as it is true for $t=0$.\\
  The system writes
 \be\label{SNL2}\bs\ba{l}\partial_tT+{\vec u}.\nabla T=uTk_0(x)\cr
 \partial_t{\vec u}+({\vec u}.\nabla){\vec u}+T\nabla Q+TQk_0(x){\vec e}_1=(1-T){\vec g}\cr
 \mbox{div}{\vec u}=0\ea\es\ee
In the system (\ref{SNL2}), appear only quadratic terms. When one wants to deduce the term of order $N$ in the system, plugging in the expansions $T^N, u^N, v^N$ and $Q^N$ one obtains source terms of the form
$$\ba{l}S_N=\sum_{j=2}^{N-1}u_jT_{N-j}k_0(x)-u_j\partial_xT_{N-j}-v_j\partial_yT_{N-j}\cr
R_{1,N}=-\sum_{j=2}^{N-1}u_j\partial_xu_{N-j}+v_j\partial_yu_{N-j}+T_j\partial_xQ_{N-j}+T_jQ_{N-j}k_0(x)\cr
R_{2,N}=-\sum_{j=2}^{N-1}u_j\partial_xv_{N-j}+v_j\partial_yv_{N-j}+T_j\partial_yQ_{N-j}\ea$$
and the system rewrites
\be\label{SNLN}\bs\ba{l}
\partial_tT_N-u_Nk_0(x)=S_N\cr
\partial_tu_N+\partial_xQ_N+Q_Nk_0(x)+gT_N=R_{1,N}\cr
\partial_tv_N+\partial_yQ_N=R_{2, N}\cr
\partial_xu_N+\partial_yv_N=0.\ea\es\ee
\paragraph{Higher order Sobolev regularity (preparatory equality)}
One of the main tools that we have to use is the divergence free condition, in order to get rid of the pressure $p$ or of the reduced pressure $Q$ when obtaining the energy inequality. Recall that the system (\ref{SNLN}) rewrites
$$\bs\ba{l}\partial_tT_N-u_Nk_0(x)=S_N\cr
\rho_0\partial_t{\vec u}_N+\nabla (\rho_0Q_N)+g\rho_0T_N{\vec e}_1=\rho_0{\vec R}_N\cr
\mbox{div}{\vec u}_N=0\ea\es$$
where ${\vec R}_N=(R_{1,N}, R_{2, N})$.\\Denote by ${\vec G}_N=\rho_0\partial_t{\vec R}_N-g\rho_0S_N{\vec e}_1$. Applying the operator $\partial_t\partial^n_{x^n}$ to equation on the velocity and using the equation on the specific volume, one obtains 
\be\label{NLeq1}\partial^n_{x^n}(\rho_0\partial^2_{t^2}{\vec u}_N)+\nabla \partial_t\partial^n_{x^n}(\rho_0Q_N)+g\partial^n_{x^n}(\rho_0k_0u_N)= \partial^n_{x^n}({\vec G}_N).\ee
One deduces the
\begin{lemma}\label{esti1}
For all $n$, one has the estimate
$$||\rho_0^{\frac12}\partial^2_{t^2}\partial^n_{x^n}{\vec u}_N||\leq C_n(\sum_{p\leq n}||\rho_0^{\frac12}\partial^p_{x^p}u_N||+ ||\rho_0^{\frac12}\partial^p_{x^p}G^N||).$$
Moreover, as the coefficients of the system depend only on $x$, this inequality is also true with the same constants when $\partial^n_{x^n}$ is replaced by $\partial^n_{x^n}\partial^q_{y^q}$ for all $q\geq 0$.
\end{lemma}
\paragraph{Proof}
One notices that (\ref{NLeq1}) writes

$$\rho_0\partial^2_{t^2}\partial^n_{x^n}{\vec u}_N+ \nabla(\partial_t\partial^n_{x^n}(\rho_0Q))+{\vec g}k_0(x)\rho_0\partial^n_{x^n}u_N= {\vec G}_N-\sum_{p=0}^{n-1}C_n^p\rho_0^{(n-p)}\partial^p_{x^p}\partial^2_{t^2}{\vec u}_N-{\vec g}\sum_{p=0}^{n-1}C_n^p\rho_0^{(n-p+1)}\partial^p_{x^p}u_N.$$
Multiplying by $\partial^2_{t^2}\partial^n_{x^n}{\vec u}_N$ and integrating, using the recurrence hypothesis that
$$||\rho_0^{\frac12}\partial^2_{t^2}\partial^p_{x^p}{\vec u}_N||\leq C_p(\sum_{m\leq p-1}||\rho_0^{\frac12}\partial^2_{t^2}\partial^m_{x^m}{\vec u}_N||)+2g^2\Lambda^2 (\sum_{m\leq p-1} ||\rho_0^{\frac12}\partial^m_{x^m}u_N||)+||G^N_{n-1}||$$
as well as the inequalities
$$|k_0(x)g|\leq \Lambda^2, |\rho_0^{-1}\rho_0^{(p)}|\leq \Lambda_p$$
(which are true as soon as $k_0$ is a $C^{\infty}$ function which derivatives are bounded, because $\rho'_0=k_0\rho_0$ )
one obtains the inequality
$$||\rho_0^{\frac12}\partial^2_{t^2}\partial^n_{x^n}{\vec u}_N||\leq C_n(\sum_{m\leq n}||\rho_0^{\frac12}\partial^2_{t^2}\partial^m_{x^m}{\vec u}_N||).$$
Lemma \ref{esti1} is proven.
\subsection{The energy equalities}
Note that the system for the leading term of the perturbation is the system (\ref{SNLN}) with a null source term. Owing to this remark, we shall treat the general case and apply the equality to the particular cases.\\
Multiplying (\ref{NLeq1}) by $\partial_t\partial^n_{x^n}{\vec u}_N$ and integrating, using the divergence free relation, one obtains
$$\ba{c}\int \partial^n_{x^n}(\rho_0\partial^2_{t^2}{\vec u}_N).\partial^n_{x^n}\partial_t{\vec u}_N dxdy+\int g\partial^n_{x^n}(\rho_0k_0u_N+\rho_0S_1^N){\vec e}_1.\partial_t\partial^n_{x^n}{\vec u}_N dxdy\cr
= \cr
\int \partial^n_{x^n}(\rho_0\partial_t{\vec S}^N).\partial_t\partial^n_{x^n}{\vec u}_N dxdy.\ea$$
In this equality, we can consider (for Sobolev inequalities) the term containing the largest number of derivatives of ${\vec u}_N$. We obtain, denoting by
$${\vec R}^N_n= \partial^n_{x^n}(\rho_0\partial^2_{t^2}{\vec u}_N)-\rho_0\partial^n_{x^n}\partial^2_{t^2}{\vec u}_N$$
$$B^N_n=\partial^n_{x^n}(\rho_0k_0u_N)-\rho_0k_0\partial^n_{x^n}u_N$$
the equality
$$\ba{c}\int\rho_0\partial^n_{x^n}\partial^2_{t^2}{\vec u}_N.\partial^n_{x^n}{\vec u}_N dxdy+\int g\rho_0k_0\partial^n_{x^n}u_N.\partial_t\partial^n_{x^n}u_N dxdy\cr
+\int{\vec R}^N_n.\partial^n_{x^n}\partial_t{\vec u}_N dxdy+\int gB_n^N.\partial_t\partial^n_{x^n}u_N dxdy\cr
= \cr
\int \partial^n_{x^n}(\rho_0\partial_t{\vec R}_N).\partial_t\partial^n_{x^n}{\vec u}_N dxdy-\int g\partial^n_{x^n}(\rho_0S_N){\vec e}_1.\partial_t\partial^n_{x^n}{\vec u}_N dxdy.\ea$$
The terms ${\vec R}_n^N$ and $B_n^N$ contain only derivatives of order less than $n-1$, hence it will appear as a source term in the application of the Duhamel principle later on. The two first terms of the previous equality are the exact derivative in time of
$$E^N_n(t)=\frac12[\int\rho_0(\partial^n_{x^n}\partial_t{\vec u}_N)^2 dxdy+\int g\rho_0k_0(\partial^n_{x^n}u_N)^2dxdy].$$
The energy equality is thus
$$E^N_n(t)=E^N_n(0)+\int_0^tg^N_n(s)ds$$
where $$\ba{ll}g_n^N(t)&=\int \partial^n_{x^n}(\rho_0\partial_t{\vec R}_N).\partial_t\partial^n_{x^n}{\vec u}_N dxdy-\int g\partial^n_{x^n}(\rho_0S_N){\vec e}_1.\partial_t\partial^n_{x^n}{\vec u}_N dxdy\cr
&=-(\int{\vec R}^N_n.\partial^n_{x^n}\partial_t{\vec u}_N dxdy+\int gB_n^N.\partial_t\partial^n_{x^n}u_N dxdy).\ea$$
Note that this source term satisfies
\be\label{NLine1}|g_n^N(t)|\leq ||\rho_0^{\frac12}\partial^n_{x^n}\partial_t{\vec u}_N||_{L^2}K_n^N(t)\ee
where one has
\be\label{NKLine1}\ba{ll}K_n^N(t)\leq &||\rho_0^{-\frac12}{\vec R}^N_n||+ ||\rho_0^{-\frac12}\partial^n_{x^n}(\rho_0\partial_t{\vec S}^N)||\cr
&+|g|[||\rho_0^{-\frac12}B_n^N||+||\partial^n_{x^n}(\rho_0S_1^N)||].\ea\ee
We are ready to prove the Duhamel inequality associated with this problem, using $gk_0(x)\leq \Lambda^2$.

\subsection{The Duhamel principle}
Two versions of the behavior of the semi group will be deduced. The first one corresponds to the general case for the terms in $\delta^2$ at least.\\
We consider the (general) system
\be\label{ss}\rho_0(x)\partial^2_{t^2}{\vec w}+\nabla(\rho_0\partial_tQ)+gk_0\rho_0w{\vec e}_1={\vec M}_0, \mbox{div}{\vec w}=0.\ee
with the initial conditions
\be\label{cis}{\vec w}\vert_{t=0}=0, \partial_t{\vec w}\vert_{t=0}=0.\ee
Note that this system is easily deduced from the system obtained for the $N$ th term of the expansion in $\delta$ of the solution.
\begin{prop}
Assume that there exists two constants $K$ and $L$, with $L>\Lambda$, sich that
\be\label{es}||\rho_0^{-\frac12}{\vec M}||\leq Ke^{Lt}.\ee
The unique solution of the linear system (\ref{ss}) with initial Cauchy conditions (\ref{cis}) satisfies the estimate
$$\bs\ba{ll}||\rho_0^{\frac12}{\vec w}||&\leq \frac{2K}{L(L-\Lambda)}(1+\frac{\Lambda^2}{(L-\Lambda)^2})^{\frac12}e^{Lt}\cr
&\leq \frac{2K}{(L-\Lambda)^2}e^{Lt}\cr
||\rho_0^{\frac12}\partial_t{\vec w}||&\leq \frac{2K}{L-\Lambda}(1+\frac{\Lambda^2}{(L-\Lambda)^2})^{\frac12}e^{Lt}\cr
||\rho_0^{\frac12}\partial^2_{t^2}{\vec w}||\leq &K(1+\frac{2\Lambda^2}{(L-\Lambda)^2})e^{Lt}\ea\es$$
\end{prop}
\paragraph{Proof}
We begin by multyplying the equation (\ref{ss}) by $\partial^2_{t^2}{\vec w}$ and integrate in space. One deduces that
$$||\rho_0^{\frac12}\partial^2_{t^2}{\vec w}||\leq \Lambda^2||\rho_0^{\frac12}w||+Ke^{Lt}.$$
We will make use of this equality later.\\
Let us multiply the equation (\ref{ss}) by $\partial_t{\vec w}$. We obtain the identity
$$\frac{d}{dt}(\frac12\int \rho_0(\partial_t{\vec w})^2dxdy+\frac12\int k_0\rho_0k_0w^2dxdy=\int M(x, y, t)\partial_t{\vec w} dxdy.$$
Integrating in time and using the initial condition (\ref{cis}) as well as the estimate (\ref{es}), we obtain the inequality
$$\int \rho_0(\partial_t{\vec w})^2dxdy\leq \Lambda^2\frac12\int k_0\rho_0w^2dxdy+2K\int_0^te^{Ls}||\rho_0^{\frac12}\partial_t{\vec w}||(s)ds.$$
Let us introduce now $u(t)=\int_0^t||\rho_0^{\frac12}\partial_t{\vec w}||(s)ds$. We obtain, considering $\frac{d}{dt}\int \rho_0{\vec w}^2dxdy$, that
$$||\rho_0^{\frac12}{\vec w}||(t)\leq \int_0^t||\rho_0^{\frac12}\partial_t{\vec w}||(s)ds$$
that is 
$$||\rho_0^{\frac12}{\vec w}||(t)\leq u(t).$$
Hence the inequality
$$(u'(t))^2\leq \Lambda^2(u(t))^2+\int_0^t2Ke^{Ls}u'(s)ds.$$
From this inequality, we deduce that $$(u'(t))^2\leq \Lambda^2(u(t))^2+2Ke^{Lt}u(t)$$
hence
$$u'(t)\leq \Lambda u(t)+\sqrt{2Ke^{Lt}u(t)}.$$
Introduce $h$ such that $u(t)=(h(t))^2e^{\Lambda t}$. We obtain the inequality
$$2hh'e^{\Lambda t}\leq \sqrt{2K}he^{\frac{L+\Lambda}{2}t}$$
hence
$$h'(t)\leq \frac12\sqrt{2K}e^{\frac{L-\Lambda}{2}t}$$
that is
$$h(t)\leq \frac{\sqrt{2K}}{L-\Lambda}e^{\frac{L-\Lambda}{2}t}$$
which leads to
$$u(t)\leq \frac{2K}{(L-\Lambda)^2}e^{Lt}$$
The estimate on $u'(t)$ follows, using $(u')^2 \leq \Lambda^2u^2+2Ke^{Lt}u$. We thus, by integration, deduce another estimate on $u$. The estimate on $\rho_0\partial^2_{t^2}{\vec w}$ is the consequence of (\ref{esti1}).\\
If one wants a general formulation of the Duhamel principle (taking into account non zero initial values), one states the following proposition, which will lead to the result of proposition \ref{improvement}, hence allowing a mixing of modes and a weak nonlinear result. The mixing of modes is not our purpose here, but we shall not speak of weak nonlinear results. See Cherfils, Garnier, Holstein \cite{garnier} for more details.\\
\begin{prop}
\label{duhamel}
The solution of
$$\frac12 \frac{d}{dt}(\int (\rho_0(\partial_t{\vec u}_N)^2-g\frac{\rho'_0}{\rho_0}\rho_0(u_N)^2)dxdy)=g(t, x, \partial_t{\vec u_N})$$
with initial condition $\partial_t{\vec u}_N(0), {\vec u_N}(0)$, with the assumption
$$|g(t, x, \partial_t{\vec u}_N)|\leq K(t)||\rho_0^{\frac12}\partial_t{\vec u}_n||_{L^2}$$
where $K$ is a positive increasing function for $t\geq 0$
satisfies the inequalities
$$\ba{l}||\rho_0^{\frac12}{\vec u}_N||^{\frac12}\leq [C_1+\int_0^t\sqrt{K(s)e^{-\Lambda s}}ds]e^{\frac{\Lambda}{2}t}\cr
||\rho_0^{\frac12}\partial_t{\vec u}_N||\leq [C_1+\int_0^t\sqrt{K(s)e^{-\Lambda s}}ds]^2e^{\Lambda t}\ea$$
where $C_1$ depends on the initial data.

\end{prop}
\paragraph{Proof}
We deduce from the energy equality the following inequality:
$$\int \rho_0(x)(\partial_t{\vec u}_N)^2dxdy-g\int k_0(x)\rho_0(x)u_N^2dxdy\leq C_{0,+}+2\int_0^tK(s)||\rho_0^{\frac12}\partial_t{\vec u}_N||_{L^2}(s)ds$$
where $C_0=\int \rho_0(x)(\partial_t{\vec u}_N)^2(0)dxdy-g\int k_0(x)\rho_0(x)u_N^2(0)dxdy$ and $C_{0, +}=\mbox{max}(C_0, 0)$.
Consider now the function  $u(t)=||\rho_0^{\frac12}{\vec u}_N(0)||+\int_0^t||\rho_0^{\frac12}\partial_t{\vec u}_N(s)||ds=||\rho_0^{\frac12}{\vec u}_N(0)||+\int_0^t||\rho_0^{\frac12}\partial_t{\vec u}_N||(s)ds$. We notice that $u'(t)=||\rho_0^{\frac12}\partial_t{\vec u}_N||(t)$ hence $u'(t)\geq 0$. Recall that $gk_0(x)\leq \Lambda^2$. The inequality implies
$$\ba{ll}(u'(t))^2&\leq \Lambda^2(u(t))^2+C_{0, +}+2\int_0^tK(s)u'(s)ds\leq\Lambda^2(u(t))^2+C_{0, +}+2K(t)u(t)\cr
&\leq (\Lambda u + \frac{K(t)}{\Lambda})^2+C_{0,+}-\frac{K(t)^2}{\Lambda^2}.\ea$$
Use now the inequality $(a^2+b^2+c^2)^{\frac12} \leq a+b+c$ for positive numbers $a, b, c$ to obtain
$$u'(t)\leq \Lambda u(t)+\sqrt{C_{0, +}}+\sqrt{2K(t)u(t)}.$$
Introducing $v(t)=u(t)e^{-\Lambda t}$ which satisfies $v(t)\geq u(0)e^{-\Lambda t}$, we deduce
$$v'(t)\leq \sqrt{C_{0, +}}e^{-\Lambda t}+\sqrt{2K(t)e^{-\Lambda t}v(t)}.$$
$\bullet$ Assume $u(0)>0$. We obtain, denoting by $h(t)=\sqrt{v(t)}$
$$2hh'\leq \sqrt{C_{0, +}}e^{-\Lambda t} +\sqrt{2K(t)e^{-\Lambda t}}h(t)$$
hence $$2h'\leq (\frac{C_{0, +}}{u(0)})^{\frac12}e^{-\frac{\Lambda t}{2}}+\sqrt{2K(t)e^{-\Lambda t}}.$$
We deduce the inequality
$$h(t)\leq h(0)+\Lambda^{-1}(\frac{C_{0, +}}{u(0)})^{\frac12}(1-e^{-\frac{\Lambda t}{2}})+\frac{1}{\sqrt2}\int_0^t\sqrt{K(s)e^{-\Lambda s}}ds.$$
which imply that there exists $A$ and $B$ such that
$$u(t)\leq (A^2e^{\Lambda t}+B^2e^{\Lambda t}(\int_0^t\sqrt{K(s)e^{-\Lambda s}}ds)^2).$$
$\bullet$ Assume $u(0)=u'(0)=0$. As $C_{0,+}=0$, we have the inequality
$$u'(t)\leq \Lambda u(t)+\sqrt{2K(t)u(t)}$$
from which one deduces, with the same notations as above, that
$$h'(t)\leq \sqrt{\frac12 K(t)e^{-\Lambda t}}$$
hence with $h(0)=0$ one obtains
$$h(t)\leq \int_0^t\sqrt{\frac12 K(s)e^{-\Lambda s}}ds.$$
$\bullet$ Assume finally $u(0)=0$ and $u'(0)>0$. We obtain
$$v'(t)\leq \sqrt{C_{0, +}}e^{-\Lambda t}+\sqrt{2K(t)e^{-\Lambda t} v(t)}.$$
Introduce ${\tilde v}(t)=v(t)-\sqrt{C_{0, +}}\frac{1-e^{-\Lambda t}}{\Lambda}$. We have
$${\tilde v}'(t)\leq \sqrt{2K(t)e^{-\Lambda t} ({\tilde v}(t)+\sqrt{C_{0, +}}\frac{1-e^{-\Lambda t}}{\Lambda})}\leq \sqrt{2K(t)e^{-\Lambda t} ({\tilde v}(t)+\frac{\sqrt{C_{0, +}}}{\Lambda})}$$
from which one deduces the inequality
$$2\sqrt{{\tilde v}+\frac{\sqrt{C_{0, +}}}{\Lambda}}\leq 2\sqrt{\frac{\sqrt{C_{0, +}}}{\Lambda}}+\int_0^t\sqrt{2K(s)e^{-\Lambda s}}ds.$$
In all the previous cases, we deduced the inequality $u(t)\leq [C_1+\int_0^t\sqrt{K(s)e^{-\Lambda s}}ds]^2e^{\Lambda t}$.\\
Using finally the relation
$$\frac{d}{dt}||\rho_0^{\frac12}{\vec u}_N||\leq ||\rho_0^{\frac 12}\partial_t{\vec u}_N||= u'(t)$$
we get
$$||\rho_0^{\frac12 }{\vec u}_N||\leq u(t)-u(0).$$
These are the two estimates of Proposition \ref{duhamel}.\\
Of course, the proof is much simpler in the case we are interested in, that is $\partial_t{\vec u}_N=0$, ${\vec u}_N=0$, where (using the notations of this paragraph, $C_0=C_{0, +}=u(0)=u'(0)=0$), where one deduces easily
$$\sqrt{u(t)e^{-\Lambda t}}\leq\int_0^t2^{-\frac12}\sqrt{K(s)e^{-\Lambda s}}ds.$$
\subsection{$H^s$ estimates for a general solution of the linearized system}
\paragraph{The $H^s$ inequalities for the solution of the homogeneous system}
We consider the system satisfied by the leading order term of the perturbation of the Euler system (which is the system (\ref{SNL2}), particular case of (\ref{SNLN}) for $N=1$. We prove in this section the analogous of the Proposition 1 of \cite{GH}, with a slightly better estimate which shows essentially that the relevant growth rate is, up to {\bf polynomial terms}, $\Lambda$:
\begin{prop}
\label{improvement}
Let $T_1(t), {\vec u}_1(t)$ be the solution of the modified linearized Euler system (\ref{SNL2}). There exists a constant $C_s$ depending only on the characteristics of the system, that is of $k_0$ and $g$, such that

$$||\rho_0^{\frac12}T_1(t)||_{H^s}+||\rho_0^{\frac12}{\vec u}_1(t)||_{H^s}\leq C_s(1+t)^s\exp(\Lambda t)(||\rho_0^{\frac12}T_1(0)||_{H^s}+||\rho_0^{\frac12}{\vec u}_1(0)||_{H^s}) .$$
\end{prop}

Note that in these inequalities (which are general) a power of $t$ appears in the bound for the norm $H^s$. This is the general case. Note that similar estimates were obtained independantly by R. Poncet \cite{poncet}.\\
An important feature of this result takes in consideration an {\bf initial} condition which is {\bf not} an eigenmode of the Rayleigh equation, and which is a combination of different eigenmodes. As we shall see in what follows, the interaction of these different eigenmodes lead to a linear growth of the form $(1+t)^se^{\Lambda t}$ for the $H^s$ norm of the solution.
\paragraph{Proof}
We prove in a first stage the $H^s$ inequality result for the system satisfied by $(T_1, u_1, v_1, Q_1)$. We use the pressure $p_1$ in the analysis.
The system imply the equation
$$\rho_0(x)\partial^2_{t^2}{\vec u}_1+\nabla\partial_tp_1=\rho_0{\vec g}k'_0u_1.$$
We apply the operator $D_{m, p}$ to this equation. The energy inequality deduced from (\ref{NLine1}) and from the inequality (\ref{NKLine1}) is
$$((u_n^1)')^2\leq \Lambda^2(u_n^1)^2+C_0+K_n^1(t)u_n^1(t)$$
where we have the estimate
$$K_n^1(t)\leq ||\rho_0^{-\frac12}{\vec R}_n^1||+|g|||\rho_0^{-\frac12}B_n^1||.$$
\begin{enumerate}
\item principal term\\
 The inequation on $||\rho_0^{\frac12}{\vec u}_1||$ writes
$$(\frac{d}{dt}||\rho_0^{\frac12}{\vec u}_1||)^2\leq \Lambda^2 ||\rho_0^{\frac12}{\vec u}_1||^2+C_0$$
hence one obtains the inequality
$$||\rho_0^{\frac12}{\vec u}_1||\leq||\rho_0^{\frac12}{\vec u}_1(0)||\cosh \Lambda t+\sqrt{\frac{C_0}{\Lambda^2}+ ||\rho_0^{\frac12}{\vec u}_1(0)||^2}\sinh\Lambda t\leq D_0e^{\Lambda t}.$$
\item derivative of the principal term\\
In the inequality obtained for $D_{1, p}{\vec u}_1$, the source term $g_1$ is bounded by $MD_0e^{\Lambda t}||D_{1, p}\partial_t{\vec u}_1||$ because it contains only derivatives of order $n-1=0$.  We have thus the inequality
$$((u_1^1)')^2\leq \Lambda^2(u_1^1)^2+C_0+2MD_0e^{\Lambda t}u_1^1(t)$$
from which one deduces
$$((u_1^1)'(t))^2\leq (\Lambda u_1^1+\frac{MD_0}{\Lambda}e^{\Lambda t})^2+C_0-(\frac{MD_0}{\Lambda})^2e^{2\Lambda t}$$
hence
$$(u_1^1)'(t)\leq \Lambda u_1^1(t)+\frac{MD_0}{\Lambda}e^{\Lambda t}+\sqrt{C_0}$$
that is
$$\frac{d}{dt}(u_1^1e^{-\Lambda t})\leq \frac{MD_0}{\Lambda}+\sqrt{C_0}e^{-\Lambda t}$$
from which one deduces
$$u_1^1(t)\leq \Lambda^{-1}(\sqrt{C_0}+MD_0t+\Lambda u_1^1(0))e^{\Lambda t}.$$
\item Greater order term:\\
We prove thus by recurrence that there exists $A_n$ and $B_n$ such that
$$u_n^1(t)\leq (A_n+B_nt)^ne^{\Lambda t}, $$
according to the inequality
$$\frac{d}{dt}(u_n^1(t)e^{-\Lambda t})\leq \sqrt{C_n}e^{-\Lambda t}+\frac{(A_{n-1}+tB_{n-1})^{n-1}}{2\Lambda}.$$
One deduces the same inequality for $\frac{d}{dt}u_n^1(t)$.
\item In the derivative $D_{n, p}$, the only term which matters for the order of the power of $t$ is $n$, hence one deduces that
$$\sum_{n+p=s}(||\rho_0^{\frac12}D_{n, p}{\vec u}_N||+||\rho_0^{\frac12}D_{n, p}\partial_t{\vec u}_N||)\leq (C_s+tD_s)^se^{\Lambda t}$$
\end{enumerate}
Proposition \ref{improvement} is proven.
Note that this improvement does not change the behavior of the approximate solution we intend to construct, because for a normal mode solution
$$u(x, y, t)={\hat u}(x)e^{iky+\gamma(k)t},$$
where $\gamma(k)$ has been calculated and where ${\hat u}(x)$ is solution of the Rayleigh equation, one has the following equalities:
\be \ba{l}||\rho_0^{\frac12}D_{m, p}{\vec u}_1(t)||= ||\rho_0^{\frac12}D_{m, p}{\vec u}_1(0)||e^{\gamma(k)t}\cr
||D_{m, p}T_1(t)||= ||T_1(0)||e^{\gamma(k)t}\cr
||\rho_0^{\frac12}D_{m, p}Q_1(t)||= ||\rho_0^{\frac12}D_{m, p}Q_1(0)||e^{\gamma(k)t}.\ea\ee
Remark that, according to Lemma \ref{existencenormalmode}, and to the equality $ikQ_1(x, y, t)=\frac{\gamma(k)}{ik}\partial_x u_1(x, y, t)$, we have also the relations
\be\ba{l}||D_{m, p}{\vec u}_1(t)||=||D_{m, p}{\vec u}_1(0)||e^{\gamma(k) t}\cr
||D_{m, p}Q_1(t)||=||D_{m, p}Q_1(0)||e^{\gamma(k) t}.\ea\ee
\subsection{The $H^s$ inequalities for the linearized system}

We consider the system (\ref{SNLN}). We apply the operator $D_{m, p}=\partial^m_x\partial^p_y$. This system becomes
\be\bs\ba{l}
\partial_tD_{m, p}T_N-D_{m, p}u_Nk_0(x)=D_{m, p}S_1^N+\sum_{q=0}^{p-1}C_p^qD_{m, q}u_Nk_0^{(q-p)}(x)\cr
\partial_tD_{m, p}u_N+\rho_0^{-1}\partial_x(\rho_0D_{m, p}Q_N)+gD_{m, p}T_N=D_{m, p}S_2^N-\sum_{q=0}^{p-1}C_p^qD_{m, q}Q_Nk_0^{(q-p)}(x) \cr
\partial_tD_{m, p}v_N+\partial_yD_{m, p}Q_N=D_{m, p}S_3^N\cr
\partial_xD_{m, p}u_N+\partial_yD_{m, p}v_N=0.\ea\es\ee
We notice that this system writes as the system (\ref{SNLN}) with a source term involving derivatives of the solution at a lesser order of derivatives in $x$.\\
We introduce
$$u_n^N(t)=||\rho_0^{\frac12}\partial_{x^n}^n{\vec u}_N(0)||_{L^2}+\int_0^t||\rho_0^{\frac12}\partial_t\partial_{x^n}^n{\vec u}_N(s)||_{L^2}ds$$
and $v_n^N(t)=u_n^N(t)e^{-\Lambda t}$.\\

We are now ready to study the behavior of the lower order terms of the expansion, assuming that we found a $\gamma(k)$ such that $\frac{\Lambda}{2}<\gamma(k)<\Lambda$.\\
We have to deal in a second part with terms of the form $u_n^N$, where $N\geq 2$. In this set-up one has to use Proposition \ref{duhamel}, because we cannot obtain the sharpest inequality using the estimate $u'\leq \Lambda u+\frac{K}{2\Lambda}+\sqrt{C_0}$.\\
\subsection{Inequalities for the following terms of the expansion}
Recall that from Lemma \ref{existencenormalmode} (proven in \cite{HelLaf}), there exists a normal mode solution of the linearized system of the form ${\hat u}(x, k)e^{iky+\gamma(k)t}$ where $\frac{\Lambda}{2}\leq \gamma(k)<\Lambda$. With this normal mode solution one constructs an approximate solution of the nonlinear system, of the form $$\ba{ll}T^N(x, y, t)=1+\sum_{j=1}^N\delta^jT_j(x, y, t)\cr
u^N(x, y, t)=\sum_{j=1}^N\delta^ju_j(x, y, t)\cr
v^N(x, y, t)=\sum_{j=1}^N\delta^jv_j(x, y, t)\cr
Q^N(x, y, t)=Q_0(x)+\sum_{j=1}^N\delta^jQ_j(x, y, t).\ea$$
There is an important Lemma, which depends on Hypothesis (H):
\begin{lemma} The functions $u_j, v_j, Q_j, T_j$ belong to $L^2$.
\end{lemma}
The proof of this Lemma is a consequence of Proposition \ref{controleQ}, which will lead to the control of the source term of the linear system on $T_N, u_N, v_N, Q_N$.\\
We shall use the estimates of Cordier, Grenier and Guo \cite{grenier}, and the method of Guo and Hwang \cite{GH} to give an $H^s$ estimate of $T^N, u^N, v^N, Q^N$ and a $L^2$ estimate of $T^N-T_0-\delta T_1$, $u^N-\delta u_1$,$v^N-\delta v_1$ to obtain a lower bound on $T^N, u^N, v^N$.\\
We prove in this section the $H^s$ estimate ${\vec u}_N$ in the  weighted norm $||\rho_0^{\frac12} .||$. Using the assumption $k_0\rho_0^{-\frac12}$ bounded, we deduce estimates in $H^s$ for ${\vec u}_N$. The first result reads as
\begin{prop}
There exists constants $C_0^p$ and $A^p$, depending only on the characteristics of the system (namely $g$, $k_0(x)$ and its derivatives) and on the $H^p$ norm of the initial data such that

$$u_p^N(t)\leq (C_0^p)^N(A^p)^{N-1}e^{N\gamma(k)t}.$$

\label{precisGH}
\end{prop}
\begin{rem}
This estimate relies heavily, as in \cite{desjardinsgrenier}, on the quadratic structure of the nonlinearity, and that we give the precise estimate on the constant $C_j$ which appears in (13) of \cite{desjardinsgrenier}. This estimate could not be obtained in the set-up of Guo and Hwang \cite{GH} because the nonlinearity was written using $\rho{\vec u}.\nabla{\vec u}$, hence a cubic nonlinearity.\end{rem}
A second comment is the following: the inequality $2\gamma(k)>\Lambda$ allows us to forget the coefficient $(1+t)^s$ in the $H^s$ estimate for a general solution of the linear system (obtained in Proposition \ref{improvement}). This is a consequence, as we shall see below, of the relation
$$e^{\Lambda t}\int_0^te^{(N\gamma(k)-\Lambda)s}ds\leq \frac{1}{N\gamma(k)-\Lambda}e^{N\gamma(k)t}$$
(to be compared with the relation $e^{\Lambda t}\int_0^te^{(\Lambda-\Lambda)s}ds\leq te^{\Lambda t}$).\\
{\bf Case $N=2$}

Recall that we have the following system
$$\bs\ba{l}\partial_tT_2-k_0(x)u_2= -u_1\partial_xT_1-v_1\partial_yT_1-T_1u_1\cr
\rho_0(x)\partial_t{\vec u}_2+\nabla (\rho_0Q_2)+gT_2=-\rho_0(x)[{\vec u}_1.\nabla {\vec u}_1]-\rho_0T_1\partial_xQ_1-\rho'_0Q_1T_1\cr
\mbox{div}{\vec u}_2=0\ea\es$$
We have thus the estimates
$$||\rho_0^{\frac12}\partial_tS_j^2||+ ||\rho_0^{\frac12}S_j^2||\leq C_j^2e^{2\gamma(k)t}.$$
This means that $K_0^2(t)\leq D_2e^{2\gamma(k)t}$, hence
$$||\rho_0^{\frac12}\partial_t{\vec u}_2||\leq (C_0+\int_0^t\sqrt{2D_2}e^{(\gamma(k)-\frac{\Lambda}{2})s}ds)^2e^{\Lambda t}$$
hence the inequality
$$||\rho_0^{\frac12}\partial_t{\vec u}_2||+||\rho_0^{\frac12}{\vec u}_2||\leq M_0e^{2\gamma(k)t}.$$
We need to derive estimates for the terms $T_2$ and $Q_2$. 
For the term $T_2$, one has
$$\frac{d}{dt}\frac12\int \rho_0T_2^2dxdy=\int k_0(x)\rho_0u_2T_2+\int S_1^2\rho_0T_2dxdy$$
from which one deduces the inequality
$$\frac{d}{dt} ||\rho_0^{\frac12}T_2||\leq M ||\rho_0^{\frac12}u_2||+ ||\rho_0^{\frac12}S_1^2||\leq(MM_0+C_1^2)e^{2\gamma(k)t}$$
hence the estimate
$$||\rho_0^{\frac12}T_2(t)||\leq ||\rho_0^{\frac12}T_2(0)||+\frac{C_1^2+MM_0}{2\gamma(k)}(e^{2\gamma(k)t}-1).$$
As for the estimate on $Q_2$, one deduces
$$\partial_x(\rho_0^{-1}\partial_x(\rho_0Q_2))+\partial^2_{y^2}Q_2+g\partial_xT_2=\mbox{div}{\vec S}^2 $$
which imply estimates on $Q_2$.\\
{\bf Case $N\geq 3$.}\\
We start with the induction hypothesis that, for $j\leq N-1$, there exists $C_0$  and $A$ such that
$$||\rho_0^{\frac12}{\vec u}_j||+||\rho_0^{\frac12}\partial_x{\vec u}_j||+||\rho_0^{\frac12}\partial_y{\vec u}_j||+|| \rho_0^{\frac12}\partial_xT_j||+|| \rho_0^{\frac12}\partial_yT_j||+|| \rho_0^{\frac12}T_j||\leq A^{j-1}C_0^je^{j\gamma(k)t}$$
and that the derivative in time of all quantities is bounded by $j\gamma(k)A^{j-1}C_0^je^{j\gamma(k)t}$. Thus there exists $M$ (independant on the number of terms which appear in the source term and which depends only on the coefficients of the system) such that the source term of (\ref{NKLine1}) for $n=0$ is bounded by:
\be\label{Ass1}K_0^N(t)\leq MA^{N-2}C_0^NN^2\gamma(k)e^{N\gamma(k)t}.\ee
Note that in this estimate the $N^2$ term comes, one from the number of the terms in the expansion $\sum_{j=0}^{N-1}A_jB_{N-j}$ and a second one\footnote{Note also that if we consider a cubic model, the number of terms in the source term is $N(N-1)$, hence adding a derivative in time we get $N^3$ in the estimate. As we can see in the following lines, this gives a less efficient estimate} from the derivative in time which appears in the source term $\partial_t{\vec S}^N$. We thus obtain, using
$$h_N(t)\leq \int_0^t\sqrt{2K_0^N(s)e^{-\Lambda s}}ds$$
the inequality
$$h_N(t)\leq\sqrt{2MA^{N-2}C_0^NN^2\gamma(k)}\int_0^t\sqrt{e^{(N\gamma(k) - \Lambda)s}}ds$$
which yields
$$h_n(t)^2\leq A^{N-1}C_0^Ne^{(N\gamma(k)-\Lambda )t}\frac{8MN^2\gamma(k)}{(N\gamma(k)-\Lambda)^2A}.$$
The choice of $A$ is thus induced by $\frac{8MN^2\gamma(k)}{(N\gamma(k)-\Lambda)^2A}\leq 1$ for all $N$ (forgetting that we have to be more precise to obtain estimates not only on ${\vec u}_N$ but also on $T_N$) hence the simplest choice is $A=\frac{8M\gamma(k)}{\gamma(k)-\frac{\Lambda}{2}}$. The value of $C_0$ is thus given by the norm of the leading term $(T_1, u_1, v_1, Q_1)$.\\
The final estimate is
$$||\rho_0^{\frac12}{\vec u}_N||\leq C_0^NA^{N-1}e^{N\gamma(k)t}$$
We proved the assumption (\ref{Ass1}).\\ 
We use this result and the estimates for a normal mode solution (on which no powers of $t$ appear for the norms of the derivatives). We obtain $$h(t)\leq h(0)+\Lambda^{-1}(\frac{C_{0, +}}{u(0)})^{\frac12}+ NC_N^{\frac12}\gamma(k)^{\frac12}\int_0^t e^{\frac{N\gamma(k)-\Lambda}{2}s}ds$$
hence as $N\gamma(k)>\Lambda$ one gets
$$h(t)\leq h(0)+\Lambda^{-1}{(\frac{C_{0, +}}{u(0)})}^{\frac12}+\frac{N}{N\gamma(k) -\Lambda}C_N^{\frac12}\gamma(k)^{\frac12}e^{\frac{N\gamma(k)-\Lambda}{2}t}.$$
We deduce the inequality (using $(a+b)^2\leq 2(a^2+b^2)$)
$$u(t)\leq 2(h(0)+u(0))+\Lambda^{-1}{(\frac{C_{0, +}}{u(0)})}^{\frac12})^2e^{\Lambda t}+2(\frac{N}{N\gamma(k) -\Lambda})^2C_N\gamma(k)e^{N\gamma(k)t}.$$
\paragraph{Remark} If the system has a cubic source term, at each stage of the construction one gets $N^{\frac12}M^{N-1}C^N$ as estimate, hence the convergence of the infinite series is not ensured by these estimates.
\subsection{Estimates for the approximate solution}
In this paragraph, we derive estimates on the global approximate solution. We shall use throughout what follows the Moser estimates, that we recall here
\be\label{M1}||D^{\alpha}(fg)||_{L^2}\leq C(||f||_\infty||g||_s+||g||_\infty||f||_s)\ee 
and
\be\label{M2}||D^{\alpha}(fg)-fD^{\alpha}g||_{L^2}\leq C(||Df||_\infty||g||_{s-1}+||g||_\infty||f||_s)\ee
and the Sobolev embedding $||f||_\infty\leq C||f||_s$ for $s>\frac{d}{2}$ and $||\nabla f||_\infty\leq C||f||_s$for $s>\frac{d}{2}+1$. More precisely, we prove that
\begin{prop}For all $\theta<1$ and for all $t<\frac{1}{\gamma(k)}\ln\frac{\theta}{\delta C_0A}$, we have
$$||T^N-1||_{H^s}+||{\vec u}^N||_{H^s}+||Q^N-q_0||_{H^s}\leq C\frac{\delta AC_0e^{\gamma(k)t}}{1-\delta AC_0e^{\gamma(k)t}}$$
$$||T^N-1||_{L^2}\geq ||T_1(0)||_{L^2}\delta e^{\gamma(k)t}-AC_0^2C_3\delta^2 \frac{e^{\gamma(k)t}}{1-\delta AC_0e^{\gamma(k)t}}$$

$$||u^N||_{L^2}\geq ||u_1(0)||_{L^2}\delta e^{\gamma(k)t}-AC_0^2C_3\delta^2 \frac{e^{\gamma(k)t}}{1-\delta AC_0e^{\gamma(k)t}}$$
$$||v^N||_{L^2}\geq ||v_1(0)||_{L^2}\delta e^{\gamma(k)t}-AC_0^2C_3\delta^2 \frac{e^{\gamma(k)t}}{1-\delta AC_0e^{\gamma(k)t}}.$$
We have also the following estimates for the remainder terms
$$||{\vec R}^N||_{H^s}+ ||S^N||_{H^s}\leq M\delta^{N+1}(N+1)^2A^{N-1}C_0^{N+2}\delta^{N+1}e^{(N+1)\gamma(k)t}.$$
\end{prop}
\paragraph{Proof}
We have proven the $H^s$ estimates for all the terms of the expansion $u_j, v_j, T_j, Q_j$. It is this easy to deduce, using (\ref{Ass1}), the estimate for the remainder terms. This comes from the inequality ($1\leq j\leq N-1$)
$$||D^{\alpha}(u_j\partial_1u_{N-j})||\leq C(||u_j||_{\infty}||u_{N-j}||_{H^{|\alpha|+1}}+ ||u_j||_{|\alpha|}||\partial_1u_{N-j}||_{\infty})$$
(and subsequent inequalities), the Solobev embedding $||f||_{\infty}\leq ||f||_2$ and the $H^s$ estimate for $s=2, 3$ for all the terms of the expansion, using also that the norm $H^s$ of the terms of the expansion in $\delta^j$ of order less than $N$ is bounded by $C_0^jA^{j-1}e^{j\gamma(k)t}$. We thus deduce that
$$||\sum_{j=1}^NT_j||_{H^s}\leq \sum_{j=1}^NCA^{j-1}C_0^j\delta^je^{j\gamma(k)t}= CC_0\delta e^{\gamma(k)t}\frac{1-(C_0A\delta)^{N-1}e^{(N-1)\gamma(k)t}}{1-C_0A\delta e^{\gamma(k)t}}.$$
When $t<T_{\delta}^{\theta}=\frac{1}{\gamma(k)}\ln\frac{\theta}{\delta C_0A}$, we obtain $1-C_0A\delta e^{\gamma(k)t}\geq 1-\theta$, hence we deduce the estimate
$$||T^N-1||_{H^s}=||\sum_{j=1}^NT_j||_{H^s}\leq \frac{CC_0}{1-\theta}\delta e^{\gamma(k)t}.$$
Moreover, one has
$$||T^N-1||_{L^2}\geq \delta ||T_1||_{L^2}-\sum_{j=2}^{N}\delta^j||T_j||$$
hence using 
$$\sum_{j=2}^{N}\delta^j||T_j||_{L^2}\leq \sum_{j=2}^{N}\delta^j||T_j||_{L^2}\leq \frac{C_0^2AC}{1-\theta}\delta^2e^{2\gamma(k)t}$$
one obtains
$$||T^N-1||_{L^2}\geq \delta ||T_1(0)||_{L^2}e^{\gamma(k)t}- \frac{C_0^2AC}{1-\theta}\delta^2e^{2\gamma(k)t}.$$
One may thus consider $C_3=\frac{C}{1-\theta}$
We thus deduce that, for $t<\frac{1}{\gamma(k)}\ln\frac{||T_1(0)||_{L^2}(1-\theta)}{CAC_0^2}$, we obtain
$$||T^N-1||_{L^2}\geq \frac12\delta ||T_1(0)||_{L^2}e^{\gamma(k)t}.$$
Similar estimates hold for $||{\vec u}^N||_{L^2}$.\\
Note that this proves that the first term of the expansion is the leading term of the approximate total solution.\\
For all what follows, we introduce
\be I(t)=\frac{A C_0e^{\gamma(k)t}}{1-\delta A C_0e^{\gamma(k)t}}\ee
\be I_{N+1}(t)=N^2A^{N-1}C_0^{N+1}e^{(N+1)\gamma(k) t}.\ee
\section{Estimates of the (nonlinear) solution}
\label{nonlinear}
We constructed in the previous section a solution $T^N, {\vec u}^N, Q^N$ such that
\be\bs\ba{ll}\partial_tT^N+{\vec u}^N\nabla T^N-k_0(x)u^NT^N= S^N\cr
\partial_t{\vec u}^N+ {\vec u}^N.\nabla{\vec u}^N+T^N\rho_0^{-1}\nabla(\rho_0Q^N)= {\vec g}+{\vec R}^N\cr
\mbox{div}{\vec u}^N=0\ea\es\ee
with the following properties for the remainder terms:
\be ||\rho_0^{\frac12}\partial^n_{x^n}R_j^N||+||\rho_0^{\frac12}\partial^n_{x^n}S^N||\leq C_n\delta^{N+1}I_{N+1}(t)\ee
\be ||\partial^n_{x^n}R_j^N||+||\partial^n_{x^n}S^N||\leq C_n\delta^{N+1}I_{N+1}(t)\ee
the constant $C_n$ depending on the Sobolev norm with weight $\rho_0^{\frac12}$ of the initial value of the normal mode solution and of the characteristic constants of the problem.\\
We deduced from this equality and the additional assumption $k_0\rho_0^{-\frac12}$ bounded that we have identical estimates on ${\vec R}^N$ and $S^N$:

\be ||\partial^n_{x^n}R_j^N||+||\partial^n_{x^n}S^N||\leq C_n\delta^{N+1}I_{N+1}(t)\ee
We study in this Section the global solution of the Euler system (\ref{SNL2}) to obtain Sobolev estimates on the difference between the approximate solution and the full solution. Let $T^d=T-T^N$, ${\vec u}^d={\vec u}-{\vec u}^N, Q^d=Q-Q^N$. We have the following system of equations:

\be\label{systemeresiduel}\bs\ba{l}\partial_tT^d+{\vec u}^N\nabla T^d+{\vec u}^d\nabla T^N= k_0(uT^d+u^dT^N)-S^N\cr
\rho_0(\partial_t{\vec u}^d+{\vec u}^d\nabla{\vec u} + {\vec u}^N\nabla{\vec u}^d)+T\nabla(\rho_0Q)-T^N\nabla(\rho_0Q^N)=-\rho_0{\vec R}^N\cr
\mbox{div}{\vec u}^d=0.\ea\es\ee

Before stating the results on the difference quantities according to the system, we use the properties of $T^N-1, {\vec u}^N, Q^N$:

\begin{lemma} 
Let $t\in [0, T^{\theta}_{\delta}]$.\\
For all $\alpha$, there exists a constant $C(|\alpha|)$ such that 
$$\ba{l}||D^{\alpha}({\vec u}^d.\nabla {\vec u}^N)||\leq C(|\alpha|)||{\vec u}^d||_{|\alpha|}\frac{CC_0\delta}{1-\theta} e^{\gamma(k)t}\cr
||D^{\alpha}(T^d(\nabla Q^N+k_0Q^N))||\leq C(|\alpha|)||T^d||_{|\alpha|}\frac{CC_0\delta}{1-\theta} e^{\gamma(k)t}\cr
||D^{\alpha}(k_0u^d(T^N-1))||\leq C(|\alpha|)||{\vec u}^d||_{|\alpha|}\frac{CC_0\delta}{1-\theta} e^{\gamma(k)t}\cr
||D^{\alpha}(T^N-1)(\nabla Q^d+k_0Q^d{\vec e}_1)||\leq C(|\alpha|)||Q^d||_{|\alpha|+1}\frac{CC_0\delta}{1-\theta} e^{\gamma(k)t}\ea$$
\end{lemma}
The proof of this Lemma comes from the fact that
$$D^{\alpha}(fg^N)=\sum C_{\alpha}^{\beta}D^{\beta}fD^{\alpha-\beta}g^N$$
and we use the estimate $||D^{\alpha-\beta}g^N||_{\infty}\leq C||g^N||_{2+|\alpha|-|\beta|}$, as well as the $H^s$ result on any term of the form $g^N=\sum_{j=1}^N \delta^jg_j$, where $g_j=u_j, v_j, T_j, Q_j$ to conclude for any term studied in the Lemma.
Moreover, we use the Moser estimates to obtain
$$||D^{\alpha}({\vec u}.\nabla f)-{\vec u}.\nabla D^{\alpha}f||\leq C(||\nabla {\vec u}||_{\infty}||\nabla f||_{|\alpha|-1}+||\nabla f||_{\infty}||{\vec u}||_{|\alpha|})$$
hence, using ${\vec u}={\vec u}^N+{\vec u}^d$, one deduces
$$||D^{\alpha}({\vec u}.\nabla f)-{\vec u}.\nabla D^{\alpha}f||\leq C(||\nabla {\vec u}^d||_{\infty}||f||_{|\alpha|}+\delta I(t)||f||_{|\alpha|}+||\nabla f||_{\infty}||{\vec u}||_{|\alpha|}+ ||\nabla f||_{\infty}\delta I(t))$$
and, similarily
$$||D^{\alpha}({\vec u}.\nabla f)||\leq C(\delta I(t)(||\nabla f||_{\infty}+||f||_{|\alpha|+1})+||{\vec u}^d||_{\infty}||f||_{|\alpha|+1}+||{\vec u}^d||_{\alpha}||\nabla f||_{\infty},$$
$$||D^{\alpha}(T\nabla Q^d)-T\nabla D^{\alpha}Q^d||\leq C(\delta I(t)||Q^d||_{|\alpha|}+
||\nabla T^d||_{\infty}||\nabla Q^d||_{|\alpha|-1}+||\nabla Q^d||_{\infty}||T^d||_{|\alpha|}$$
according to the equality $D^{\alpha}(T\nabla Q^d)-T\nabla D^{\alpha}Q^d=D^{\alpha}((T-1)\nabla Q^d)-(T-1)\nabla D^{\alpha}Q^d$.

We shall also use the following estimates
\be\label{unonlineaire}||D^{\alpha}({\vec u}^d.\nabla {\vec u}^d)||\leq C(|\alpha|)||{\vec u}^d||_4||{\vec u}^d||_{|\alpha|+1}\ee
\be\label{unonlineaire2}||D^{\alpha}({\vec u}^d.\nabla {\vec u}^d)-{\vec u}^d.\nabla D^{\alpha}{\vec u}^d||\leq C(|\alpha|)||{\vec u}^d||_4||{\vec u}^d||_{|\alpha|}.\ee
These equalities come respectively from (\ref{M1}) and (\ref{M2}).

Introduce in what follows ${\vec V}={\vec u}^d.\nabla{\vec u}^N+{\vec u}^N.\nabla{\vec u}^d$, ${\vec W}={\vec V}+T^d\rho_0^{-1}\nabla (\rho_0Q^N)$. We have the estimates
$$\ba{l}||D^{\alpha}{\vec V}||\leq M_{|\alpha|}I(t)\delta ||{\vec u}^d||_{|\alpha|+1}\cr
||D^{\alpha}{\vec W}||\leq M_{|\alpha|}I(t)\delta (||{\vec u}^d||_{|\alpha|+1}+ ||T^d||_{|\alpha|}), \forall \alpha\cr
||D^{\alpha}(T^d\nabla Q^d)-T^d\nabla D^{\alpha}Q^d||\leq C(||T^d||_3||Q^d||_2+||T^d||_4||Q^d||_1)\mbox{ for }|\alpha|=2\cr
||T^d\rho_0^{-\frac12}\nabla(\rho_0Q^N)||\leq \delta I(t)||T^d||\ea$$
\subsection{Estimates on the density}
The equation on the density yield
$$\partial_tT^d+{\vec u}.\nabla T^d-k_0uT^d=k_0u^dT^N-{\vec u}^d.\nabla T^N-S^N.$$
Apply the operator $D^{\alpha}$ and denote by $W^1_{\alpha}=D^{\alpha}({\vec u}.\nabla T^d)-{\vec u}.\nabla D^{\alpha}T^d$. This equation
rewrites
$$\partial_tD^{\alpha}T^d+{\vec u}.\nabla D^{\alpha}T^d+W^1_\alpha-D^{\alpha}(k_0uT^d)+D^{\alpha}({\vec u}^d.\nabla
T^N)-D^{\alpha}(k_0u^dT^N)=0.$$
We can decompose $W^1_\alpha-D^{\alpha}(k_0uT^d)$ into two parts, the one with ${\vec u}^N$, the other one with ${\vec u}^d$, denoted
respectively by $W_\alpha$ and $W^N_\alpha$. It is clear that
$$||W^N_\alpha-D^{\alpha}(k_0uT^d)+D^{\alpha}({\vec u}^d.\nabla
T^N)-D^{\alpha}k_0u^dT^N)||\leq C\delta (||{\vec u}_d||_{|\alpha|}+||T^d||_{|\alpha|}).$$
it is also clear that, using Moser estimates, $||W_\alpha||\leq C(||\nabla {\vec u}^d||_\infty||T^d||_{|\alpha|}+||{\vec
u}^d||_{|\alpha|}||T^d||_{\infty})$. One is thus left with the inequality

$$\frac{d}{dt}||D^{\alpha}T^d||\leq ||W_\alpha||+||D^{\alpha}S^N||+C\delta I(t)(||{\vec
u}_d||_{|\alpha|}+||T^d||_{|\alpha|}).$$
We have thus the estimate
$$\ba{ll}\frac{d}{dt}||D^{\alpha}T^d||&\leq C(||\nabla {\vec u}^d||_\infty||T^d||_{|\alpha|}+||{\vec
u}^d||_{|\alpha|}||T^d||_{\infty})+\delta^{N+1}I_{N+1}M\cr
&+C\delta I(t)(||{\vec
u}_d||_{|\alpha|}+||T^d||_{|\alpha|}).\ea$$

\subsection{Estimates on the pressure}
We obtained the relations
$$||\nabla Q^d||\leq M_1(||{\vec u}^d.\nabla {\vec u}^d||+\delta I(t)(||{\vec u}^d||_1+||T^d||)+ \delta^{N+1}I_{N+1}(t)).$$
$$\ba{ll}\sum_{|\alpha=1}||\nabla D^{\alpha}Q^d||&\leq M_2(\sum_{|\alpha|=1}||D^{\alpha}({\vec u}^d.\nabla {\vec u}^d)||+ \delta
I(t)[||{\vec u}^d||_2+||T^d||_1]\cr
&+ \delta^{N+1}I_{N+1}(t)
+(1+\delta I(t)+||T^d||_3)(||{\vec u}^d.\nabla {\vec u}^d||\cr
&+\delta I(t)(||T^d||+||{\vec
u}^d||_1)+\delta^{N+1}I_{N+1}(t))\ea$$
Using the fact that $t\leq T^\delta$, one obtains

$$||\nabla Q^d||\leq M_1(||{\vec u}^d.\nabla {\vec u}^d||+||{\vec u}^d||_1+||T^d||+ \delta^{N+1}I_{N+1}(t)).$$
$$\ba{ll}\sum_{|\alpha=1}||\nabla D^{\alpha}Q^d||\leq &M_2(\sum_{|\alpha|=1}||D^{\alpha}({\vec u}^d.\nabla {\vec u}^d)||+ 
||{\vec u}^d||_2+||T^d||_1\cr
&+(1+||T^d||_3)(||{\vec u}^d.\nabla {\vec u}^d||+||T^d||+||{\vec
u}^d||_1)+\delta^{N+1}I_{N+1}(t))\ea$$

In what follows, we introduce
$${\vec G}^{\alpha}_N=D^{\alpha}(T^N\nabla Q^d)-T^N\nabla D^{\alpha}Q^d+D^{\alpha}(T^NQ^dk_0{\vec e}_1)+ D^{\alpha}(T^d\nabla
Q^N+k_0T^dQ^N{\vec e}_1)+ D^{\alpha}{\vec R}^N,$$
$${\vec G}^{\alpha}=D^{\alpha}(T^d\nabla Q^d)-T^d\nabla D^{\alpha}Q^d+D^{\alpha}(T^dQ^dk_0{\vec e}_1).$$
The equation on $D^{\alpha}{\vec u}^d$ is
$$\partial_tD^{\alpha}{\vec u}^d+D^{\alpha}({\vec u}^d.\nabla {\vec u}^d)+ D^\alpha {\vec R}^N+{\vec G}^{\alpha}+{\vec
G}^{\alpha}_N+T\nabla D^{\alpha}Q^d=0.$$
When one multiplies by $\nabla D^{\alpha}Q^d$, one uses the divergence free condition on $D^{\alpha}{\vec u}^d$ to get the estimate
$$\frac23||\nabla D^{\alpha}Q^d||\leq ||D^{\alpha}({\vec u}^d.\nabla {\vec u}^d)||+ ||D^{\alpha}{\vec R}^N||+ ||{\vec G}^{\alpha}||+||{\vec
G}^{\alpha}_N||.$$
We use $$||{\vec
G}^{\alpha}_N||\leq C(1+t)^{|\alpha|+3}(||Q^d||_{|\alpha|}+||T^d||_{|\alpha|}+||{\vec u}^d||_{|\alpha|+1})\delta I(t)$$
and
$$||{\vec
G}^{\alpha}||\leq C(||\nabla T^d||_{\infty}||Q^d||_{|\alpha|}+||T^d||_{|\alpha|}||\nabla
Q^d||_{\infty}+||T^d||_{\infty}||Q^d||_{|\alpha|}).$$
Hence we obtain (and it is pertinent for $|\alpha|>2$)
$$\ba{ll}||\nabla D^{\alpha}Q^d||\leq &C'(||D^{\alpha}({\vec u}^d.\nabla {\vec u}^d)||+ ||D^{\alpha}{\vec R}^N||)+ C\delta
I(t)(||Q^d||_{|\alpha|}+||T^d||_{|\alpha|}+||{\vec u}^d||_{|\alpha|+1})\cr
&+ C(||T^d||_3||Q^d||_{|\alpha|}+||T^d||_{|\alpha|}||Q^d||_3).\ea$$
For $|\alpha|=2$, we will obtain $||Q^d||_3$, which is important.

We use the equality, for $|\alpha|=2$
$$||{\vec G}^{\alpha}||=||D^{\alpha}T^d\nabla Q^d+\sum_{0<\beta<\alpha}D^{\beta}T^d\nabla D^{\alpha-\beta}Q^dC^{\beta}_{\alpha}||$$
which leads to the inequality
$$||{\vec G}^{\alpha}||\leq D_0(||T^d||_4||Q^d||_1+||T^d||_3||Q^d||_2).$$
Replacing this estimate in the inequality for $\alpha$ such that $|\alpha|=2$, one gets
$$||D^{\alpha}\nabla Q^d||\leq C_1(||D^{\alpha}({\vec u}^d.\nabla {\vec u}^d)||+||T^d||_2+||{\vec u}^d||_3+
(1+||T^d||_3)||Q^d||_2+||T^d||_4||Q^d||_1+ \delta^{N+1}I_{N+1}(t).$$
Using the inequalities on $||Q^d||_1$ and $||Q^d||_2$, one gets
$$||Q^d||_1\leq M_1(||{\vec u}^d.\nabla {\vec u}^d||+||{\vec u}^d||_1+||T^d||+ \delta^{N+1}I_{N+1}(t))$$
$$||Q^d||_2\leq M_2(||{\vec u}^d.\nabla {\vec u}^d||_1+||{\vec u}^d||_2+||T^d||_1+
\delta^{N+1}I_{N+1}(t)(1+||T^d||_3)+(1+||T^d||_3)(1+||T^d||+||{\vec u}^d.\nabla {\vec u}^d||))$$
$$\ba{ll}||Q^d||_3&\leq M_3(||{\vec u}^d.\nabla {\vec u}^d||_2+||{\vec u}^d||_3+||T^d||_2+ 
(1+||T^d||_3^2+||T^d||_4)||{\vec u}^d.\nabla {\vec u}^d||+ (1+||T^d||_3)||{\vec u}^d||_2\cr
&+||T^d||_4||{\vec u}^d||+
\delta^{N+1}I_{N+1}(t)(1+||T^d||_4+(1+||T^d||_3)^2))\ea$$

We use then the inequalities
$$\ba{ll}||\nabla D^{\alpha}Q^d||&\leq C(||D^{\alpha}({\vec u}^d.\nabla {\vec u}^d)||+||D^{\alpha}{\vec R}^N||+ ||Q^d||_{|\alpha|}\cr
&+||T^d||_{|\alpha|}(1+||Q^d||_3)+||{\vec u}^d||_{|\alpha|+1}+(1+||T^d||_3)||Q^d||_{|\alpha|})\ea$$
from which one obtains
$$\ba{ll}||Q^d||_{|\alpha|+1}&\leq M_{|\alpha|+1}(||{\vec u}^d.\nabla {\vec u}^d||_{|\alpha|}+||{\vec u}^d||_{|\alpha|+1}+\delta^{N+1}I_{N+1}(t)\cr
&+||T^d||_{|\alpha|}(1+||Q^d||_3)+||Q^d||_{|\alpha|}(1+||T^d||_3)).\ea$$
Note that we have the estimate
\be ||{\vec u}^d.\nabla {\vec u}^d||_{|\alpha|}\leq C||{\vec u}^d||_3||{\vec u}^d||_{|\alpha|+1}.\ee
hence
\be \ba{ll}||Q^d||_{|\alpha|+1}&\leq M_{|\alpha|+1}((1+||{\vec u}^d||_3)||{\vec u}^d||_{|\alpha|+1}+\delta^{N+1}I_{N+1}(t)\cr
&+||T^d||_{|\alpha|}(1+||Q^d||_3)+(1+||T^d||_3)||Q^d||_{|\alpha|}\ea\ee
It is then enough to use a recurrence argument to control the norm of $Q^d$ in $H^{s+1}$ using the control of the norm of $Q^d$ in $H^s$.

For the control on ${\vec u}^d$, let us rewrite the equation on $D^{\alpha}{\vec u}^d$. We introduce 
$${\vec V}_{\alpha}=D^{\alpha}({\vec u}.\nabla {\vec u}^d)-{\vec u}.\nabla D^{\alpha}{\vec u}^d, {\vec W}_{\alpha}=D^{\alpha}(T.\nabla Q^d)-T.\nabla D^{\alpha}Q^d
.$$
We have the estimates
$$||{\vec V}_{\alpha}||\leq C(1+||{\vec u}^d||_3)||{\vec u}^d||_{|\alpha|}, ||{\vec W}_{\alpha}||\leq C(1+||T^d||_3)||Q^d||_{|\alpha|}.$$
Using the relation
$$\int T\nabla D^{\alpha}Q^dD^{\alpha}{\vec u}^ddxdy=-\int D^{\alpha}Q^d(\nabla(T^N-1)+\nabla T^d)D^{\alpha}{\vec u}^ddxdy$$
thanks to the divergence free condition, as well as
$$\int {\vec u}.\nabla D^{\alpha}{\vec u}^d.D^{\alpha}{\vec u}^ddxdy=0$$
one obtains the estimate
$$\frac{d}{dt}||D^{\alpha}{\vec u}^d||\leq ||{\vec V}_{\alpha}||+||{\vec W}_{\alpha}||+ ||D^{\alpha}{\vec R}^N||+ ||D^{\alpha}(k_0TQ^d)||+||D^{\alpha}Q^d||(1+||T^d||_3),$$
hence the inequality
\be\frac{d}{dt}||D^{\alpha}{\vec u}^d||\leq C[(1+||T^d||_3)||Q^d||_{\alpha}+(1+||{\vec u}^d||_3)||{\vec u}^d||_{|\alpha|}+\delta^{N+1}I_{N+1}(t)]\ee

For $|\alpha|\geq 3$, this inequality is an {\it a priori} inequality. We have to state the identical inequalities for $|\alpha|=0, 1, 2$.\\
We have the following inequalities:
\be\frac{d}{dt}||{\vec u}^d||\leq C_0((1+||T^d||_3)||Q^d||+\delta^{N+1}I_{N+1}(t))\ee
because ${\vec V}_{\alpha}={\vec W}_{\alpha}=0$,
\be\frac{d}{dt}||\nabla {\vec u}^d||\leq C((1+||T^d||_3)||Q^d||_1+(1+||{\vec u}^d||_3)||{\vec u}_d||_1+(1+||T^d||_1)||Q^d||_3+\delta^{N+1}I_{N+1}(t))\ee
and
\be\frac{d}{dt}||{\vec u}^d||_2\leq C((1+||T^d||_3)||Q^d||_2+(1+||T^d||_4)||Q^d||_1+(1+||{\vec u}^d||_3)||{\vec u}^d||_2+ \delta^{N+1}I_{N+1}(t))\ee

We thus deduce an estimate of the form

$$\frac{d}{dt}(||T^d||^2_4+||{\vec u}^d||^2_4)\leq C(1+||T^d||_4)^4+||{\vec u}^d||_3)(||T^d||^2_4+||{\vec u}^d||^2_4)+ \delta^{N+1}I_{N+1}(t)(||T^d||^2_4+||{\vec u}^d||^2_4)^{\frac12}$$
from which one deduces an estimate of the form
$$\frac{d}{dt}(||T^d||^2_4+||{\vec u}^d||^2_4)^{\frac12}\leq C(1+||T^d||_4)^4+||{\vec u}^d||_3)(||T^d||^2_4+||{\vec u}^d||^2_4)^{\frac12}+C\delta^{N+1}I_{N+1}(t)$$

\paragraph{End of the proof}
We thus know that, for $t\leq T^{\delta}$, we have $\delta^{N+1}I_{N+1}(t)\leq 1$ hence an inequality of the form
$$\frac{d}{dt}H(t)\leq C((1+(H(t))^4)H(t)+1)$$
where $H(t)=(||T^d||^2_4+||{\vec u}^d||^2_4)^{\frac12}$.\\
As we have $H(0)=0$, one deduces that
$$\int_0^{H(t)}\frac{ds}{(1+s^4)s+1}\leq Ct.$$
The function $H\rightarrow \int_0^{H}\frac{ds}{(1+s^4)s+1}$ is a bijection from $[0, +
\infty[$ onto $[0, \int_0^{+\infty}\frac{ds}{(1+s^4)s+1}[$. For $H(t)\geq 1$, one deduces $Ct\geq  \int_0^{1}\frac{ds}{(1+s^4)s+1}$, hence for $t<\frac{1}{C}\int_0^{1}\frac{ds}{(1+s^4)s+1}=T_1$, one obtains $H(t)\leq 1$. The set of points $t$ such that $t>0$ and $H(t)\leq 1$ is not empty.\\
Once this set is not empty (and once we proved that the solution exists for a time $T_1$), we obtain
\begin{lemma}
\label{lemmefinal}
Let $h$ be a function such that
$$\frac{dh}{dt}\leq C(1+h(t))^4h(t)+C\delta^{N+1}e^{(N+1)\gamma(k)t}, h(0)=0.$$
For $\delta<1$ and $(N+1)\gamma(k)>17C$, denoting by $T_0^{\delta}=\frac{1}{\gamma(k)}\ln\frac{1}{\delta}$, one has

$$\forall t\in [0, T_0^{\delta}], h(t)\leq \delta^{N+1}e^{(N+1)\gamma(k)t}.$$
\end{lemma}
\paragraph{Proof}
The inequality we start with is
$$\frac{d}{dt}h(t)\leq C(1+h(t))^4h(t)+C\delta^{N+1}e^{(N+1)\gamma(k)t}.$$
We consider $N$ such that $(N+1)\gamma(k)>17C$. We study the interval where $h(t)\in [0, 1]$, knowing that $h(0)=0$. Consider
$t_0$ the first time (if it exists) where $h(t_0)=1$. If it does not exist, then
$h(t)\leq 1$ for $t\in [0, T^{\delta}_0]$ and we have, for all $t\in [0, T^{\delta}_0]$
 the inequality
 
 $$h'(t)\leq 16Ch(t)+C\delta^{N+1}e^{(N+1)\gamma(k)t}.$$
 from which one deduces
 $$h(t)\leq \frac{C\delta^{N+1}}{(N+1)\lambda -16C}e^{(N+1)\gamma(k)t}<\delta^{N+1}e^{(N+1)\gamma(k)t}$$
 hence $h(T_0^{\delta})<1$.\\
 If $t_0$ exists, we have, for all $t\in [0, t_0]$, the inequality
 $$\frac{d}{dt}(h(t)e^{-16Ct})\leq -C(1-h(t))h(t)R(h(t))e^{-16Ct}+
 C\delta^{N+1}e^{(N+1)\gamma(k) t-16Ct}$$
 where $R(x)=(1+x)^3+2(1+x)^2+4(1+x)+8$, from which one deduces that
 $$h(t_0)e^{-16Ct_0}\leq
 \frac{C}{(N+1)\gamma(k)-16C}\delta^{N+1}e^{(N+1)\gamma(k)
 t_0-16Ct_0}<\delta^{N+1}e^{(N+1)\gamma(k)
 t_0-16Ct_0}$$
 hence $h(t_0)<1$, contradiction.\\
 We thus deduce that $h(t)\leq 1$ for $t\in [0, T^{\delta}_0]$, hence
 
 $$h(t)\leq \delta^{N+1}e^{(N+1)\gamma(k)t}, t\in [0, T^{\delta}_0].$$
 Lemma \ref{lemmefinal} is proven.\\
 We have thus the inequalities
 $$||{\vec u}||\geq ||{\vec u}^N||-||{\vec u}^d||\geq \delta ||{\vec u}_1(0)||-C_0^2A\delta^2\frac{e^{2\gamma(k)t}}{1-C_0\delta Ae^{\gamma(k)t}}-\delta^{N+1}e^{(N+1)\gamma(k)t}.$$
 Choose $t=T_1^{\delta}=\frac{1}{\gamma(k)}\ln\frac{\theta }{C_0A\delta}$. We have
 $$||{\vec u}||\geq \delta e^{\gamma(k)t}[||{\vec u}_1(0)||-C_0\frac{\theta}{1-\theta}-\theta^N].$$
 We thus check that there exists $\varepsilon_0\leq \frac56$ such that $\theta<\varepsilon_0$ implies $[||{\vec u}_1(0)||-C_0\frac{\theta}{1-\theta}-\theta^N]\geq\frac12||{\vec u}_1(0)||$. Hence for $t\leq \frac{1}{\gamma(k)}\ln\frac{\varepsilon_0 }{C_0A\delta}$, one has
 \be ||{\vec u}(t)||\geq \frac12||{\vec u}_1(0)||\delta e^{\gamma(k)t}.\ee
In particular
 $$||{\vec u}(\frac{1}{\gamma(k)}\ln\frac{\varepsilon_0 }{C_0A\delta})||\geq \frac{\varepsilon_0}{2}||{\vec u}_1(0)||.$$
 We proved Theorem \ref{resNL}.\\
It is then clear that, for $0\leq t\leq T^{\delta}_1$, this term is smaller than $\theta$,as small as one wants, hence the inequality on $T^d, {\vec u}^d$.\\
As $T=T^N+T^d, {\vec u}={\vec u}^N+{\vec u}^d$, one obtains
$$||T-1||_{\infty}\geq ||T^N||-||T^d||$$
which imply the result.

\end{document}